\documentclass[smallextended]{svjour3}
\smartqed  
\usepackage{amsfonts}
\usepackage{graphicx}
\usepackage{amssymb,amsxtra,amsmath}
\usepackage{latexsym}
\usepackage{booktabs}
\usepackage{cite}
\usepackage{url}
\usepackage{color}
\usepackage[section]{placeins}

\newcommand{\bicgstab}{{BICGSTAB}}

\newcommand{\calP}{{ \mathcal P}}
\newcommand{\calI}{{ \mathcal I}}
\def\calG{{ \mathcal G}}

\newcommand{\vX}{ \mathbf{x} }
\newcommand{\vY}{ \mathbf{y} }
\newcommand{\vg}{ \mathbf{g} }

\newcommand{\lf}{\left}
\newcommand{\rt}{\right}
\newcommand{\lp}{\left(}
\newcommand{\rp}{\right)}
\newcommand{\ep}{\varepsilon}

\newcommand{\unrml}{\mathbf{n}}
\newcommand{\calp}{\mathbf{p}}

\begin{document}

\title{A Radial Basis Function (RBF)-Finite Difference (FD) Method for Diffusion and Reaction-Diffusion Equations on Surfaces}
\titlerunning{An RBF-FD Method for Reaction-Diffusion Equations on Surfaces}

\author{Varun Shankar\and Grady B. Wright\and Robert M. Kirby\and Aaron L. Fogelson}
\authorrunning{V. Shankar, G.B. Wright, R.M. Kirby and A.L. Fogelson}

\institute{Varun Shankar \at
              School of Computing, University of Utah, Salt Lake City, UT 84112 \\
              \email{shankar@cs.utah.edu}           
           \and
           Grady B. Wright \at
              Department of Mathematics, Boise State University, Boise, ID 83725-1555 \\
              \email{gradywright@boisestate.edu}
           \and
           Robert M. Kirby \at
              School of Computing, University of Utah, Salt Lake City, UT 84112 \\
              \email{kirby@sci.utah.edu}           
           \and   
           Aaron L. Fogelson \at
           		Department of Mathematics, University of Utah, Salt Lake City, UT 84112 \\
           		\email{fogelson@math.utah.edu}
}
\date{Received: date / Accepted: date}
\maketitle
\begin{abstract}
In this paper, we present a method based on Radial Basis Function (RBF)-generated Finite Differences (FD) for numerically solving diffusion and reaction-diffusion equations (PDEs) on closed surfaces embedded in $\mathbb{R}^d$. Our method uses a method-of-lines formulation, in which surface derivatives that appear in the PDEs are approximated \emph{locally} using RBF interpolation. The method requires only scattered nodes representing the surface and normal vectors at those scattered nodes. All computations use only extrinsic coordinates, thereby avoiding coordinate distortions and singularities. We also present an optimization procedure that allows for the stabilization of the discrete differential operators generated by our RBF-FD method by selecting shape parameters for each stencil that correspond to a global target condition number. We show the convergence of our method on two surfaces for different stencil sizes, and present applications to nonlinear PDEs simulated both on implicit/parametric surfaces and more general surfaces represented by point clouds.
\end{abstract}

\section{Introduction}
\label{sec:intro}

Methods based on global Radial Basis Functions (RBFs) have become quite popular for the numerical solution of the partial differential equations (PDEs) due to their ability to handle scattered node layouts, their simplicity of implementation and their spectral accuracy and convergence on smooth problems. While these methods have been successfully applied to the solution of PDEs on planar regions~\cite{Fasshauer:2007}, they have also been applied to PDEs on the two-sphere $\mathbb{S}^2$ (\emph{e.g.}~\cite{Gia:2005, FlyerWright:2007, FlyerWright:2009}).

Many methods have been developed for the solution of the class of PDEs known as diffusion (or reaction-diffusion) equations on more general surfaces. Of these, the so-called \emph{intrinsic methods} attempt to solve PDEs using surface-based meshes and coordinates intrinsic to the surface under consideration; this approach can be efficient since the dimension of the discretization is restricted to the dimension of the surface under consideration (\emph{e.g.}~\cite{Calhoun2009,Dziuk2007}). However, such intrinsic coordinates can contain singularities or distortions which are difficult to accomodate. A popular alernative is the class of so-called \emph{embedded, narrow-band methods} that extend the PDE to the embedding space, construct differential operators in extrinsic coordinates, and then restrict them to a narrow band around the surface (\emph{e.g.}~\cite{MacDonaldRuuth2009,Piret2012}). Such methods incur the additional expense of solving equations in the dimension of the embedding space; the curse of dimensionality will ensure these costs will grow rapidly depending on the order of accuracy of the method. 

RBFs have recently been used to compute an approximation to the surface Laplacian in the context of a pseudospectral method for reaction-diffusion equations on manifolds~\cite{FuselierWright2012}. In that study, global RBF interpolants were used to approximate the surface Laplacian at a set of ``scattered'' nodes on a given surface, combining the advantages of intrinsic methods with those of the embedded methods. This method showed very high rates of convergence on smooth problems on parametrically and implicitly defined manifolds. However, for $N$ points on the surface, the cost of that method scales as $O(N^3)$. Furthermore, the dense nature of the resulting differentiation matrices means that the cost of applying those matrices to solution vectors is $O(N^2)$, assuming the manifold is static. Our goal is to develop a method that is less costly to apply than the global RBF method while still retaining the ability to use scattered nodes on the surface to approximate derivatives, thereby combining the benefits of the intrinsic and narrow-band approaches. Our motivation is to eventually apply this method for the simulation of chemical reactions on evolving surfaces of platelets and red blood cells. For this, we turn to RBF-generated Finite Differences (RBF-FD).

First discussed by Tolstykh~\cite{Tolstykh2003}, RBF-FD formulas are generated from RBF interpolation over \emph{local} sets of nodes on the surface.  This type of method is conceptually similar to the standard FD method with the exception that the differentiation weights enforce the exact reproduction of derivatives of shifts of RBFs (rather than derivatives of polynomials as is the case with the
standard FD method) on each local set of nodes being considered. This results in sparse matrices like in the standard FD method, but with the added advantage that the RBF-FD method can naturally handle
irregular geometries and scattered node layouts. We note that the RBF-FD method has proven successful for a number of other applications in planar domains in two and higher dimensions (\emph{e.g.}~\cite{ShuDing2003,CecilQian2004,Wright200699,Chandhini2007,SPLM}). The RBF-FD method has also been shown to be successful on the surface of a sphere~\cite{FoL11,FlyerLehto2012} for convective flows by stabilization with hyperviscosity. 

An RBF-FD method for the solution of diffusion and reaction-diffusion equations on general 1D surfaces embedded in 2D domains was recently developed~\cite{ShankarWrightEtAlIJNMF2013}. In our experiments, a straightforward extension of that approach to 2D surfaces proved to be unstable, requiring hyperviscosity-based stabilization as in the case of the RBF-FD method for purely convective flows. In this work, we modify the RBF-FD formulation presented in \cite{ShankarWrightEtAlIJNMF2013}, and present numerical and algorithmic strategies for generating RBF-FD operators on general surfaces. Our approach appears to do away with the need for hyperviscosity-based stabilization. 

The remainder of the paper is organized as follows. In Section 2, we briefly review RBF interpolation of both scalar and vector data on scattered node sets in $\mathbb{R}^d$. Section 3 discusses the formulation of surface differential operators in Cartesian coordinates. Section 4 then goes on to describe how these differential operators are discretized in the form of \emph{sparse} differentiation matrices and presents a method-of-lines formulation for the solution of diffusion and reaction-diffusion equations on surfaces; this section also presents important implementation details and comments on the computational complexity of our RBF-FD method. In Section 5, we detail our shape parameter optimization approach and illustrate how it can be used to stabilize the RBF-FD discretization of the surface Laplacian without the need for hyperviscosity-based stabilization. In Section 6, we numerically demonstrate the convergence of our method for different stencil sizes (on two different surfaces) for the forced scalar diffusion equation using two different approaches to selecting the shape parameter $\ep$. Section 7 demonstrates applications of the method to simulations of Turing Patterns on two classes of surfaces: implicit and parametric surfaces, and more general surfaces represented only by point clouds. We conclude our paper with a summary and discussion of future research directions in Section 8.

\textbf{Note}: Throughout this paper, we will use the terms surface or manifold to refer to smooth embedded submanifolds of codimension one in $\mathbb{R}^d$ with no boundary, with the specific case of $d = 3$. Although not pursued here, straightforward extensions are possible for manifolds of higher codimension, or manifolds of codimension $1$ embedded in higher or lower dimensional spaces.

\section{A review of RBF interpolation}
\label{sec:rbf_review}

We start with a review of RBF interpolation, which is essential to understanding the RBF-FD approach outlined in the next section. Let $\Omega \subseteq \mathbb{R}^d$, and $\phi :\Omega \times\Omega \to \mathbb{R}$ be a kernel with the property $\phi(\vX,\vY) := \phi(\|\vX-\vY\|)$ for $\vX,\vY\in\Omega$, where $\|\cdot\|$ is the standard Euclidean norm in $\mathbb{R}^d$.  We refer to kernels with this property as \emph{radial kernels} or \emph{radial functions}. Given a set of nodes $X = \{\vX_k\}_{k = 1}^N \subset \Omega$ and a continuous target function $f:\Omega \to \mathbb{R}$ sampled at the nodes in $X$, we consider constructing an RBF interpolant to the data of the following form:
\begin{align}
I_{\phi} f(\vX) = \sum_{k=1}^N c_k \phi(\|\vX - \vX_k\|) + c_{N+1}.
\label{eq:srbfinterp}
\end{align}
The interpolation coefficients $\{c_k\}_{k=1}^{N+1}$ are determined by enforcing $\lf.I_{\phi} f\rt|_{X} = \lf.f\rt|_{X}$ and $\sum_{k=1}^{N} c_k = 0$. This can be expressed as the following linear system:
\begin{align}
\underbrace{
\begin{bmatrix}
\phi(r_{1,1}) & \phi(r_{1,2}) & \hdots & \phi(r_{1,N}) & 1 \\
\phi(r_{2,1}) & \phi(r_{2,2}) & \hdots & \phi(r_{2,N}) & 1\\
\vdots & \vdots & \ddots & \vdots & \vdots\\
\phi(r_{N,1}) & \phi(r_{N,2}) & \hdots & \phi(r_{N,N}) & 1\\
1 & 1 & \hdots & 1 & 0
\end{bmatrix}}_{A_X}
\underbrace{
\begin{bmatrix}
c_1 \\
c_2 \\
\vdots \\
c_N \\
c_{N+1}
\end{bmatrix}}_{c_f}
=
\underbrace{
\begin{bmatrix}
f_1 \\
f_2 \\
\vdots \\
f_N\\
0
\end{bmatrix}}_{f_X},
\label{eq:rbf_linsys}
\end{align}
where $r_{i,j} = ||\vX_i - \vX_j||$. If $\phi$ is a positive-definite radial kernel or an order one conditionally positive-definite kernel on $\mathbb{R}^d$, and all nodes in $X$ are distinct, then the matrix $A_{X}$ above is guaranteed to be invertible (see, for example,~\cite[Ch.\ 6--8]{Wendland:2004}). 

In the present study, we are interested in the set of interpolation nodes $X$ lying on a lower dimensional surface $\Omega=\mathbb{M}$ in $\mathbb{R}^d$.  However, we will still use the standard Euclidean distance in $\mathbb{R}^d$  for $\| \cdot \|$ in Equation \eqref{eq:srbfinterp} (\emph{i.e.,} straight line distances rather than distances intrinsic to the surface).  This significantly simplifies constructing interpolants as no explicit information about the surface is needed.  A theoretical foundation for RBF interpolation on surfaces with this distance measure is given in~\cite{FuselierWright:2010}, where the authors prove and demonstrate that favorable error estimates can be achieved.

In describing our method for approximating the surface Laplacian in the next section, it is useful to extend the above discussion to the interpolation of vector-valued functions $\vg(\vX): \Omega \to \mathbb{R}^d$ sampled at a set of nodes $X = \{\vX_k\}_{k = 1}^N \subset \Omega$. For this problem, we simply apply scalar RBF interpolation as given in Equation \eqref{eq:srbfinterp} to each component of $\vg(\vX)$ and represent the resulting interpolant as $I_{\Phi}\vg$.  For example, if $d=3$ and $\vg = \begin{bmatrix} g^x & g^y & g^z \end{bmatrix}^T$, then the vector interpolant is given as
\begin{align}
I_{\Phi} \vg(\vX) = \begin{bmatrix} I_{\phi} g^x(\vX) & I_{\phi} g^y(\vX) & I_{\phi} g^z(\vX) \end{bmatrix}.
\label{eq:vrbfinterp}
\end{align}
The interpolation coefficients for each component of $I_{\Phi} \vg$ can be determined by solving a system of equations similar to the one listed in Equation \eqref{eq:rbf_linsys}, but with the right-hand-side replaced with the respective component of $\vg$ sampled on $X$.  This allows some computational savings for determining the interpolation coefficients for $I_{\phi} f$ and $I_{\Phi} \vg$ with a direct solver since the matrix $A_{X}$ then only needs to be factored once.

There are many choices of positive definite or order one conditionally positive definite radial kernels that can be used in applications; see~\cite[Ch. 4, 8, 11]{Fasshauer:2007} for several examples. These kernels can be classified into two types: finitely smooth and infinitely smooth. It is still an open question as to which kernel is optimal for which application.  Typically infinitely smooth kernels such as the Gaussian ($\phi(r) = \exp(-(\ep r)^2))$, multiquadric ($\phi(r) = \sqrt{1 + (\ep r)^2}$), and inverse multiquadric ($\phi(r) = 1/\sqrt{1 + (\ep r)^2}$) are used in the RBF-FD method for numerically solving PDEs~\cite{ShuDing2003,Wright200699,Bayona2010,Davydov2011,FlyerLehto2012}.  We continue with this trend in the present work and use the inverse multiquadric (IMQ) kernel, which is positive definite in $\mathbb{R}^d$, for any $d$.  

All infinitely smooth kernels, features a free ``shape parameter'' $\ep$, which can be used to change the kernels from peaked (large $\ep$) to flat (small $\ep$). In the limit as $\ep\to 0$ (i.e. a flat kernel), RBF interpolants to data scattered in $\mathbb{R}^d$ typically (and always in the case of the Gaussian radial kernel) converge to (multivariate) polynomial interpolants~\cite{DriscollFornberg2002,LarFor05,Sch05}, and, in the case of the surface of a sphere, they converge to spherical harmonic interpolants~\cite{FoPi07}.  For smooth target functions, smaller (but non-zero) values of $\ep$ generally lead to more accurate RBF interpolants~\cite{FoWr,LarFor05}.  However, the standard way of computing these interpolants by means of solving Equation \eqref{eq:rbf_linsys} (referred to as RBF-Direct in the literature) becomes ill-conditioned for small $\ep$ (see, \emph{e.g.,}~\cite{FZ07}).  While some stable algorithms have been developed for bypassing this ill-conditioning~\cite{FoWr,FoPi07,FaMC12,FLF,FoLePo13}, there are issues with applying them to problems where the interpolation nodes are arranged on a lower dimensional surface than the embedding space, as is the case in the present study.  These issues are related to the nodes being ``non-unisolvent'' and some strategies have recently been undertaken to resolve them~\cite{LLHF}, but a robust approach is not yet available.  In later sections of this study, we will detail strategies for selecting $\ep$ based on condition numbers of RBF interpolation matrices. We will also introduce a strategy for modifying $\ep$ to produce interpolants that compensate for irregularities in point spacing on our test surfaces.  


\section{Surface Laplacian in Cartesian coordinates}
Here we review how to express the surface Laplacian in Cartesian (or extrinsic) coordinates; for a full discussion see~\cite{FuselierWright2012}.  Working with the operator in Cartesian coordinates is fundamental to our proposed method as it completely avoids singularities that are associated with using intrinsic, surface-based coordinates ({\em e.g.} the pole singularity in spherical coordinates).  We restrict our discussion to surfaces $\mathbb{M}$ of dimension two embedded in $\mathbb{R}^3$ since these are the most common in applications.

Let $\calP$ denote the projection operator that takes an arbitrary vector field in $\mathbb{R}^3$ at a point
$\vX=(x,y,z)$ on the surface and projects it onto the tangent
plane to the surface at $\vX$.  Letting $\unrml=(n^x,n^y, n^z)$ denote
the \emph{unit} normal vector to the surface at $\vX$, this
operator is given by
\begin{align}
\calP = \calI - \unrml\unrml^{T} =
\begin{bmatrix}
(1 - n^x n^x) & -n^x n^y & -n^x n^z  \\
-n^x n^y & (1 - n^y n^y) & -n^y n^z \\
-n^x n^z & -n^y n^z & (1-n^z n^z) 
\end{bmatrix} =
\begin{bmatrix}
\calp^x & 
\calp^y &
\calp^z 
\end{bmatrix},
\label{eq:proj_mat}
\end{align}
where $\calI$ is the $3$-by-$3$ identity matrix, and $\calp^x$, $\calp^y$ and $\calp^z$
are vectors representing the projection operators in the $x$, $y$ and $z$
directions, respectively.  We can combine $\calP$ with the standard
gradient operator in $\mathbb{R}^3$, $\nabla
= \begin{bmatrix} \partial_x & \partial_y & \partial_z \end{bmatrix}^T$, to define
the \emph{surface} gradient operator $\nabla_{\mathbb{M}}$ in Cartesian
coordinates as
\begin{align}
\nabla_{\mathbb{M}} := \calP\nabla = \begin{bmatrix} \calp^x \cdot \nabla \\ \calp^y \cdot \nabla \\ \calp^z \cdot \nabla \end{bmatrix} = \begin{bmatrix} \calG^x \\ \calG^y \\ \calG^z \end{bmatrix}.
\label{eq:grad}
\end{align}
Noting that the surface Laplacian $\Delta_{\mathbb{M}}$ is given as the surface divergence of the surface gradient, this operator can be written in Cartesian coordinates as
\begin{align}
\Delta_{\mathbb{M}}:= \nabla_{\mathbb{M}}\cdot \nabla_{\mathbb{M}} = (\calP\nabla) \cdot \calP\nabla = \calG^x \calG^x + \calG^y \calG^y + \calG^z \calG^z.
\label{eq:Lap}
\end{align}
The approach we use to approximate the surface Laplacian mimics the
formulation given in Equation \eqref{eq:Lap} and is conceptually similar to
the approach based on global RBF interpolation used in~\cite{FuselierWright2012}, with the important difference being that we use local RBF interpolants.

\section{RBF-FD approximation to the surface Laplacian}
\label{sec:MoL}

Let $X = \{\vX_k\}_{k=1}^N$ denote a set of (scattered) node locations on a surface $\mathbb{M}$ of dimension two embedded in $\mathbb{R}^3$ and suppose $f:\mathbb{M}\rightarrow\mathbb{R}$ is some differentiable function sampled on $X$. Our goal is to approximate $\lf.\Delta_{\mathbb{M}} f\rt|_X$ with finite-difference-style local approximations to the operator $\Delta_{\mathbb{M}}$. Without loss of generality,
let the node where we want to approximate $\Delta_{\mathbb{M}}f$ be
$\vX_1$, and let $\vX_2,\ldots,\vX_n$ be the $n-1$ nearest
neighbors to $\vX_1$, measured by Euclidean distance in $\mathbb{R}^3$. We refer to $\vX_1$ and its $n-1$ nearest neighbors as the \emph{stencil} on the surface corresponding to $\vX_1$ and denote this stencil as $P_1 = \{\vX_k\}_{k=1}^n$.  We seek an approximation to $\Delta_{\mathbb{M}}f$ at $\vX_1$ that involves a linear combination of the values of $f$ over the stencil $P_1$ of the form
\begin{align} 
(\Delta_{\mathbb{M}}f)\bigr|_{\vX=\vX_1} \approx \sum_{j=1}^{n} w_j f(\vX_j).
\label{eq:fd_weights}
\end{align}
The weights $\{w_j\}_{j=1}^n$ in this approximation will be computed using RBFs, and will be referred to as RBF-FD weights.

The first step to computing the RBF-FD weights is to construct an RBF interpolant of $f$ similar to Equation \eqref{eq:srbfinterp}, but now only over the nodes in $P_1$, \emph{i.e.}
\begin{align}
I_{\phi}f(\vX) = \sum_{j=1}^n c_j\phi(r_j(\vX)) + c_{n+1},
\label{eq:f_interp}
\end{align}
where $r_j(\vX) = \|\vX-\vX_j\|$.  The interpolation coefficients $c_j$ can be determined by the solution to the system of equations given in Equation \eqref{eq:rbf_linsys}, but with $X$ replaced with $P_1$;  we denote this system by $A_{P_1} c_f = f_{P_1}$.   Second, we compute the surface gradient of the above interpolant using Equation \eqref{eq:grad} and evaluate it at the nodes in $P_1$.  In the case of the $\calG^x$ component of the gradient, this is given as
\begin{align}
\lf( \calG^x I_{\phi}f(\vX) \rt)\bigr|_{\vX = \vX_i} = \sum_{j=1}^n c_j \underbrace{\lf(\calG^x \phi(r_j(\vX)) \rt)\bigr|_{\vX=\vX_i}}_{{\lf(B^x_{P_1}\rt)}_{i,j}}\; , i=1,\ldots,n,
\label{eq:grad_x_f_interp}
\end{align}
where the constant term from Equation \eqref{eq:f_interp} has vanished since its gradient is zero.  We can rewrite Equation \eqref{eq:grad_x_f_interp} in matrix-vector form using the fact that $c_f = A_{P_1}^{-1} f_{P_1}$ as follows:
\begin{align}
\lf( \calG^x I_{\phi}f \rt)\bigr|_{P_1} = B^x_{P_1} c_f = \lf(B^x_{P_1} A^{-1}_{P_1}\rt)f_{P_1} = G_{P_1}^x f_{P_1}.
\end{align}
Here $G_{P_1}^x$ is an $n$-by-$n$ differentiation matrix that represents the RBF approximation to the $x$-component of the surface gradient operator over the set of nodes in $P_1$.  Similar approximations can be obtained to the $y$- and $z$-components of the surface gradient operator on this stencil as follows:
\begin{align}
\lf( \calG^y I_{\phi}f \rt)\bigr|_{P_1} = \lf(B^y_{P_1} A^{-1}_{P_1}\rt)f_{P_1} = G_{P_1}^y f_{P_1}, \\
\lf( \calG^z I_{\phi}f \rt)\bigr|_{P_1} = \lf(B^z_{P_1} A^{-1}_{P_1}\rt)f_{P_1} = G_{P_1}^z f_{P_1},
\end{align}
where the entries of $B^y_{P_1}$ and $B^z_{P_1}$ are given as
\begin{align*}
(B^y_{P_1})_{i,j} = \lf(\calG^y \phi(r_j(\vX))\rt)\bigr|_{\vX=\vX_i} \; \text{and}\; (B^z_{P_1})_{i,j} = \lf(\calG^z \phi(r_j(\vX))\rt)\bigr|_{\vX=\vX_i}.
\end{align*}
In the third step, we mimic the continuous formulation of the surface Laplacian in Equation \eqref{eq:Lap} using the differentiation matrices $G_{P_1}^x$, $G_{P_1}^y$, and $G_{P_1}^z$ in place of the operators $\calG^x$, $\calG^y$, and $\calG^z$, respectively, which gives the following approximation to the surface Laplacian of $f$ at all the nodes in $P_1$:
\begin{align}
\lf(\Delta_{\mathbb{M}} f\rt)\bigr|_{P_1} \approx \underbrace{\lf(G_{P_1}^x G_{P_1}^x + G_{P_1}^y G_{P_1}^y + G_{P_1}^z G_{P_1}^z \rt)}_{L_{P_1}}f_{P_1}.
\label{eq:surface_lap_P1}
\end{align}
This approximation is equivalent to the following operations: construct an interpolant of $f$ over $P_1$, compute its surface gradient, interpolate each component of the surface gradient, apply the surface divergence, and evaluate it at $P_1$.  Hence, we can use the vector interpolant notation from Equation \eqref{eq:vrbfinterp}, to write Equation \eqref{eq:surface_lap_P1} equivalently as
\begin{align*}
\lf(\Delta_{\mathbb{M}} f\rt)\bigr|_{P_1} \approx \lf(\nabla_{\mathbb{M}}\cdot I_{\Phi} \lf(\nabla_{\mathbb{M}}I_{\phi}f\rt) \rt) \bigr|_{P_1}.
\end{align*}
This approach of repeated interpolation and differentiation avoids the need to analytically differentiate the surface normal vectors of $\mathbb{M}$, which implies closed form expressions for these values are not needed. This simplifies the computations and makes the method applicable to surfaces defined by point clouds (as illustrated in Section \ref{sec:results2}).

While the approximation in Equation \eqref{eq:surface_lap_P1} is for all the nodes in $P_1$,  we are only interested in the approximation at $\vX=\vX_1$ (the ``center'' point of the stencil $P_1$) according to Equation \eqref{eq:fd_weights}.  Because of the ordering of nodes in $P_1$, the value of Equation \eqref{eq:fd_weights} is given by the first value in the vector that results from the product on the right of Equation \eqref{eq:surface_lap_P1}.  Thus, the weights $w_j$ in Equation \ref{eq:fd_weights} are given by the entries in the first row of the matrix $L_{P_1}$ from Equation \eqref{eq:surface_lap_P1}.  Extracting these entries from this matrix, and disregarding the rest, then completes the steps for determining the RBF-FD weights for the node $\vX_1$.  

For each node $\vX_j\in X$, $j=1,\ldots,N$, we repeat the above procedure of finding its $n-1$ nearest neighbors (stencil $P_j$), computing the corresponding matrix $L_{P_j}$ according to Equation \eqref{eq:surface_lap_P1}, and extracting out of this matrix the row of RBF-FD weights for $\vX_j$.  These weights are then arranged into a \emph{sparse} $N$-by-$N$ differentiation matrix $L_X$ for approximating the surface Laplacian over all the nodes in $X$.  

The computational cost of computing each matrix $L_{P_j}$ is $O(n^3)$, and there are $N$ such stencils, so that the total cost of computing the entries of $L_X$ is $O(n^3 N)$ (this is apart from the cost of determining the stencil nodes, for which an efficient method is discussed below).   In practice, $n << N$ and would typically be fixed as $N$ increases, so that the total cost scales like $O(N)$. Furthermore, each $L_{P_j}$ can be computed independently from the others and is thus a embarrassingly parallel computation.  In contrast, the method from \cite{FuselierWright2012}, requires $O(N^3)$ operations and results in a dense differentiation matrix.  However, the accuracy of this global method is better than the local RBF-FD approach.

\subsection{Implementation details}
\label{sec:impl}
To efficiently determine the members of stencils $P_k$, $k=1,\ldots,N$, we first build a k-d tree for the full set of $N$ nodes in $X$. The k-d tree is constructed in $O(d N \log N)$ operations, where $d$ is the number of dimensions. The members of stencil $P_k$ can then be determined from the k-d tree in $O(\log N)$ operations.  Combining the computational cost of the k-d tree construction and look-ups with computing the RBF-FD weights, the total cost of building $L_X$ is $O(N \log N) + O(N)$, where the constants in the last term depend on the cube of $n$.

\begin{figure}[ht]	
\centering
\includegraphics[width=2.3in,height=1.3in]{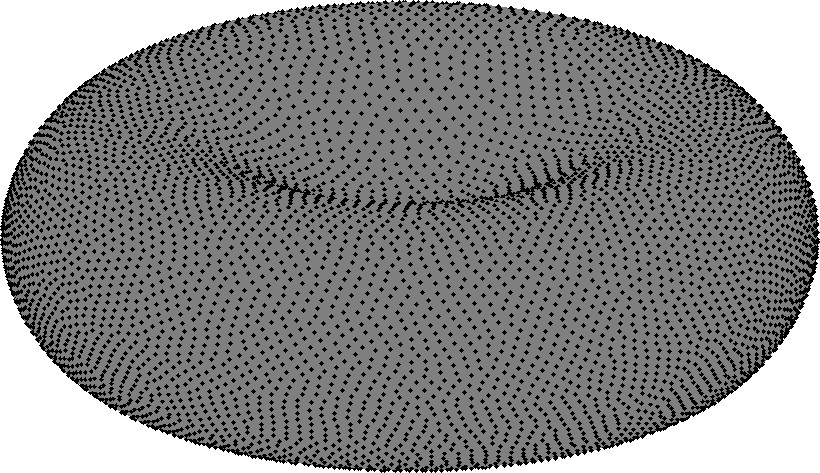}
\includegraphics[width=2.3in,height=2.0in]{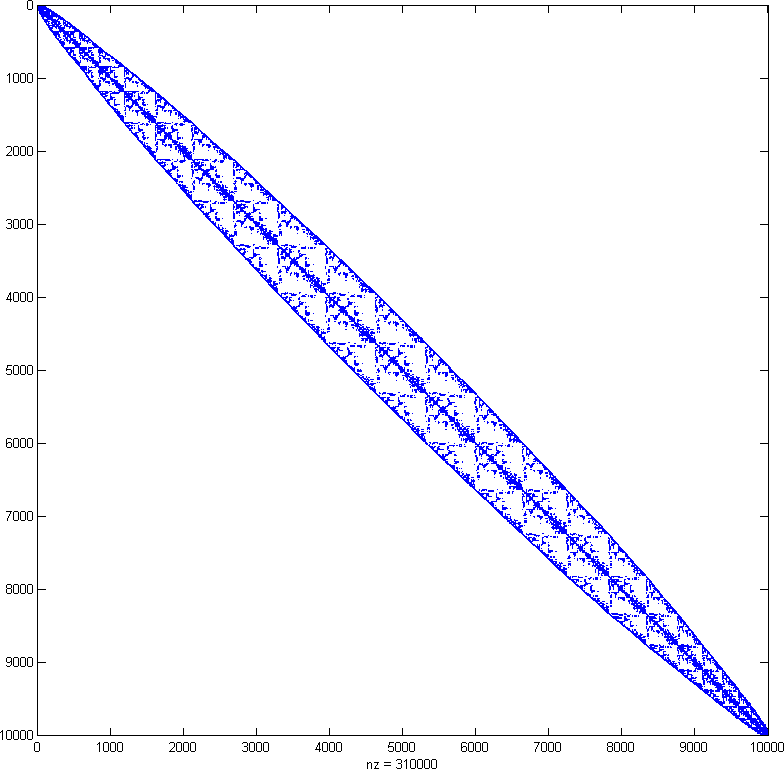}

\vspace{0.4cm}

\includegraphics[width=2.3in,height=2.2in]{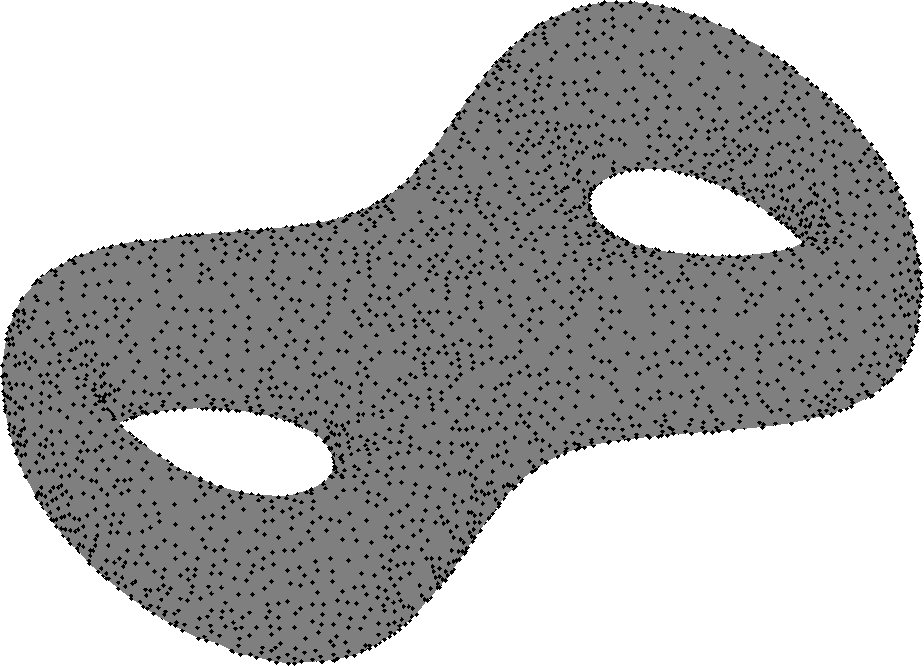}
\includegraphics[width=2.3in,height=2.0in]{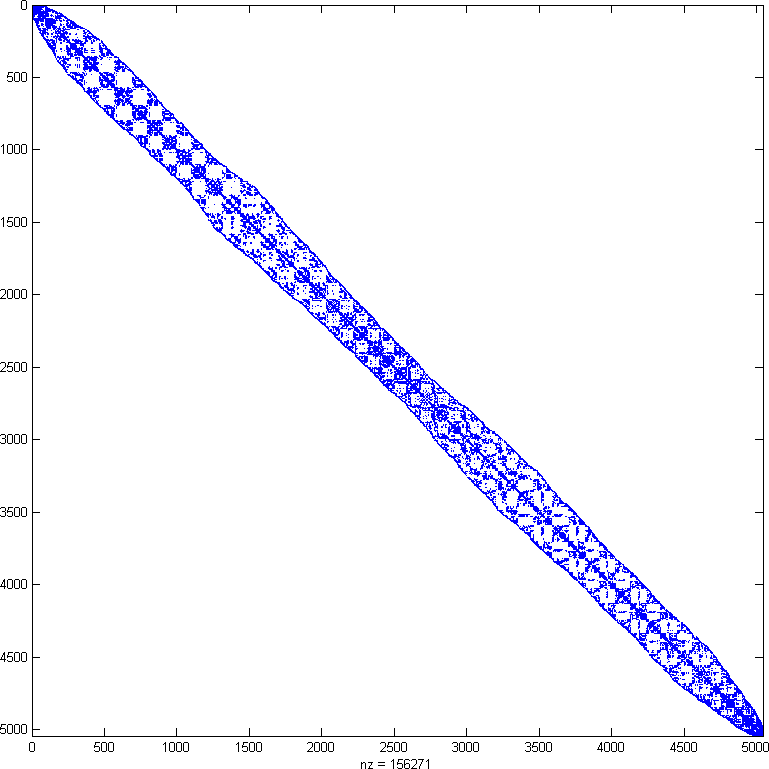}
\caption{The figure on the top left shows maximum determinant nodes mapped from the sphere to the surface of the Red Blood Cell. The figure on the top right shows the re-ordered matrix $L_X$, obtained by applying the Reverse Cuthill-McKee re-ordering algorithm. $0.31\%$ of the entries of the matrix are non-zeros for the Red Blood Cell. The figure on the bottom left shows a node set obtained on the double-torus. The figure on the bottom right shows the re-ordered matrix $L_X$, obtained by applying the Reverse Cuthill-McKee re-ordering algorithm. $0.62\%$ of the entries of the matrix are non-zeros for the double-torus. We use a stencil size of $n = 31$ for both objects.}
\label{fig:rcm}	
\end{figure}

 We note that points in each of the $N$ stencils $P_k$ on the surfaces are selected merely using a distance criterion; in other words, for a node $\vX_k$, the stencil only comprises of its $n-1$ nearest neighbors, with the distances measured in $\mathbb{R}^3$, rather than along the surface. While it is possible to include more information to form more regular or biased stencils, we do not explore these possibilities in our current work. One consequence of using distances in the embedding space is that one must exercise caution when simulating PDEs on the surfaces of thin objects, or thin features of more general surfaces. If the distance between points across a thin feature is smaller than the distance between points on the same sides of the surface of the thin feature, a poor approximation to the surface Laplacian will result in that region. We will not address this issue in our work, except by taking care to have a sufficiently dense sampling of the surface around thin features.

In general, the nodes in $X$ will lack any ordering, which may negatively impact the fill-in of the sparse matrix $L_X$. We therefore re-order the matrix using the Reverse Cuthill-McKee algorithm~\cite{GeorgeLiuRCM} for all our tests. This usually results in faster iterations in an iterative solver involving $L_X$, and improved sparsity in stored lower triangular and upper triangular factors within a sparse direct solver involving $L_X$.  Figure \ref{fig:rcm} shows two different surfaces, a idealized red blood cell and a double-torus, with corresponding nodes $X$, and the resulting sparsity pattern of $L_X$ after re-ordering with Reverse Cuthill-McKee.

\subsection{Method-of-lines}

In Sections \ref{sec:convergence} and \ref{sec:results2}, we use the RBF-FD discrete approximation to the surface Laplacian in the method-of-lines (MOL) to simulate diffusion and reaction-diffusion equations on surfaces.  We briefly review this technique for the former equation, as its generalization to the latter follows naturally.

The diffusion of a scalar quantity $u$ on a surface with a (non-linear) forcing term is given as
\begin{align}
\frac{\partial u}{\partial t} = \nu \Delta_{\mathbb{M}} u + f(t,u),
\label{eq:diffusion_scalar}
\end{align}
where $\delta>0$ is the diffusion coefficient, $f(t,u)$ is the forcing term, and an initial value of $u$ at time $t=0$ is given.  Letting $X=\{\vX_j\}_{j=1}^N\subset\mathbb{M}$ and $u_X \in\mathbb{R}^N$ denote the vector containing the samples of $u$ at the points in $X$, our RBF-FD method for \eqref{eq:diffusion_scalar} takes the form
\begin{align}
\frac{d}{dt} u_X = \delta L_X u_X + f\lp t,u_X\rp,
\label{eq:kernel_mol}
\end{align}
where $L_X$ is an $n$-node RBF-FD differentiation matrix for approximating $\Delta_{\mathbb{M}}$ over the nodes in $X$, as described above.   This is a (sparse) system of $N$ coupled ODEs and, provided it is stable (see Section \ref{sec:stability}), can be advanced in time with a suitably chosen time-integration method.  For an explicit time-integration method, $L_X$ can be evaluated in $O(N)$ operations.  For a method that treats the diffusion term implicitly, one can use an iterative solver such a $\bicgstab$, or form the sparse upper and lower triangular factors obtained from the LU factorization of the implicit equations and use them for an efficient direct solver every time-step.  These are the two respective approaches we use in our convergence studies in Section \ref{sec:convergence} and our applications in Section \ref{sec:results2}.

We conclude by noting that solving surface reaction-diffusion equations with an RBF-FD method was also considered in our paper~\cite{ShankarWrightEtAlIJNMF2013} for the case of 1D surfaces embedded in $\mathbb{R}^2$.  However, the approach used in that study for computing a discrete approximation to the surface Laplacian differs in an important way from the RBF-FD formulation of the surface Laplacian presented above.  In that work, given a set of $N$ nodes ($X$) on a surface, we start by using $n$-node RBF-FD formulas to construct differentiation matrices for the $\calG^x$ and $\calG^y$ over the node set $X$, which we denote by $G^x_X$ and $G^y_X$.  Next, the surface Laplacian was approximated from these matrices as $L_X =  G^x_X G^x_X + G^y_X G^y_X$. As with the above approach, this formulation also avoids the need to compute derivatives of the normal vectors of the surfaces, but has the effect of doubling the bandwidth of the $L_X$ compared to $G^x_X$ and $G^y_X$.  We tried extending this approach to two dimensional surfaces in embedded in $\mathbb{R}^3$, but encountered stability issues when combining this with the method-of-lines, as the differentiation matrices $L_X$ had eigenvalues with (sometimes large) positive real parts. The present method appears to be much less susceptible to these problems as discussed in the next section.


\section{Shape Parameter and Eigenvalue Stability}
\label{sec:stability}
A necessary condition for stability of the MOL approach described in the previous section is that the eigenvalues of the RBF-FD differentiation matrices $L_X$ must be in the stability domain of the ODE solver used for advancing the system in time.  As a minimum requirement, this will generally mean that all eigenvalues must, at the very least, be in the left half plane. The RBF-FD procedure does not guarantee that this property will hold for $L_X$, and it is possible to encounter situations in which this requirement is violated. In this section, we discuss a procedure related to choosing a stencil-dependent shape parameter $\ep_k$ when computing the RBF-FD weights that appears to ameliorate this issue and lead to $L_X$ with eigenvalues in the left-half plane. 

\begin{figure}[ht]	
\centering
\includegraphics[width=2.3in,height=2.2in]{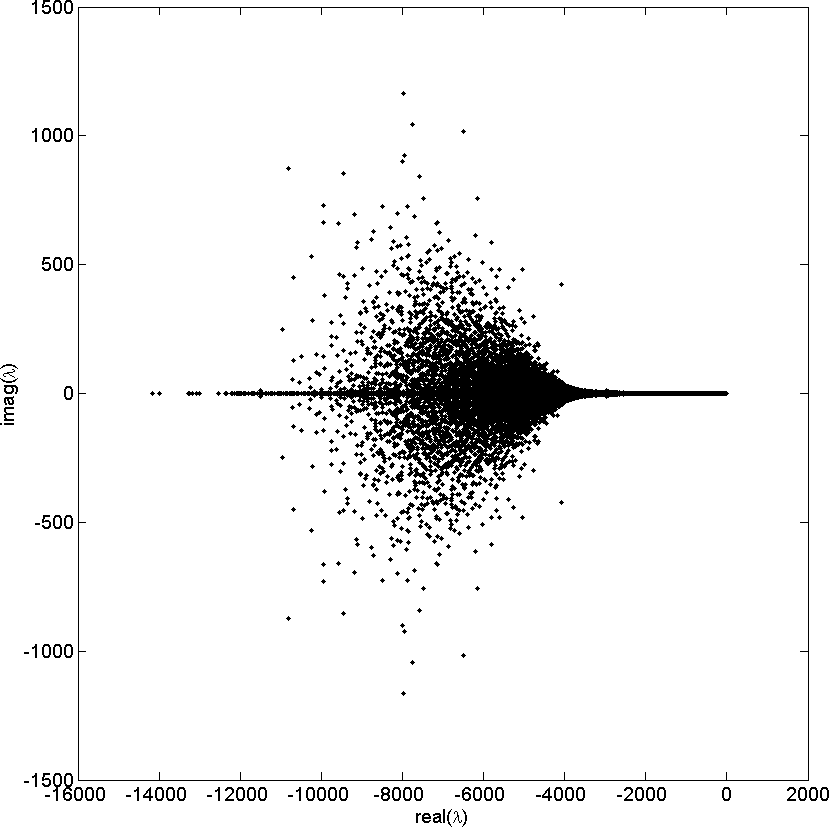}
\includegraphics[width=2.3in,height=2.2in]{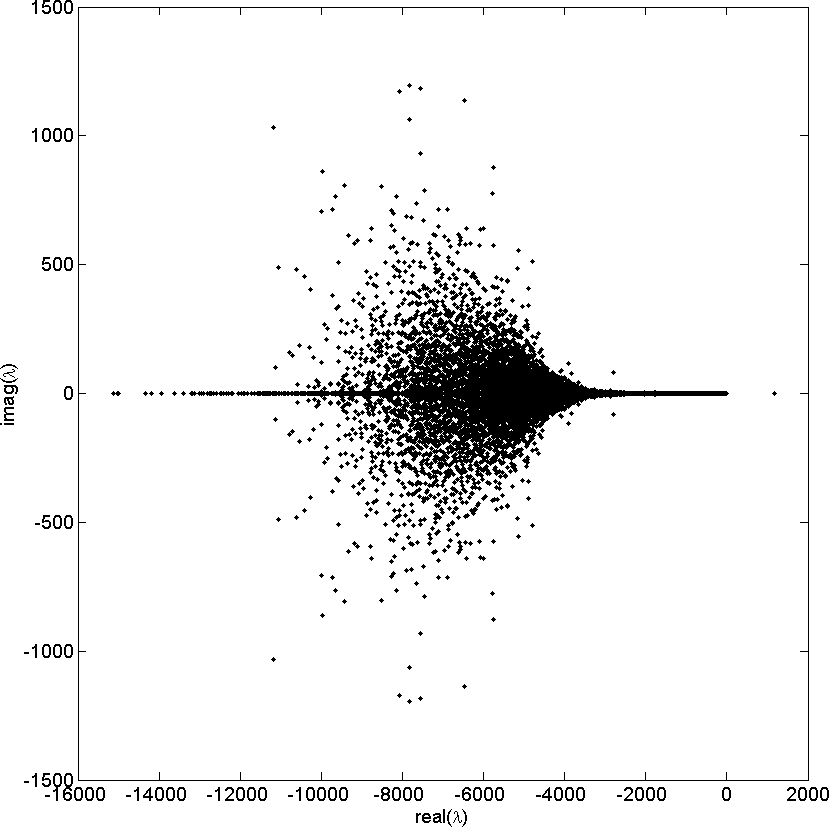}

\vspace{0.4cm}

\includegraphics[width=2.3in,height=2.2in]{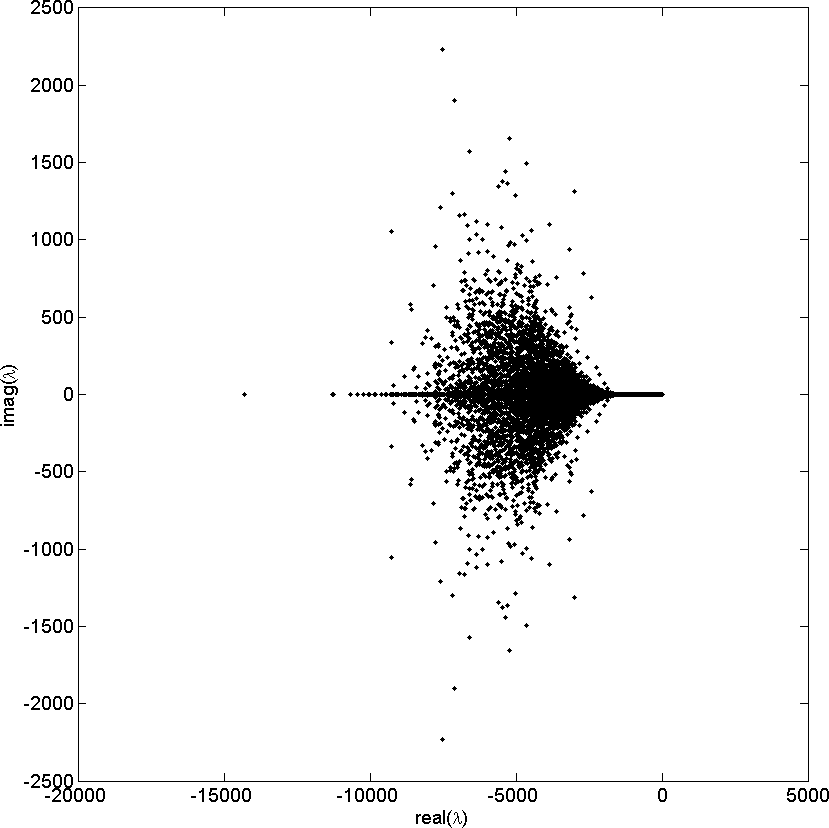}
\includegraphics[width=2.3in,height=2.2in]{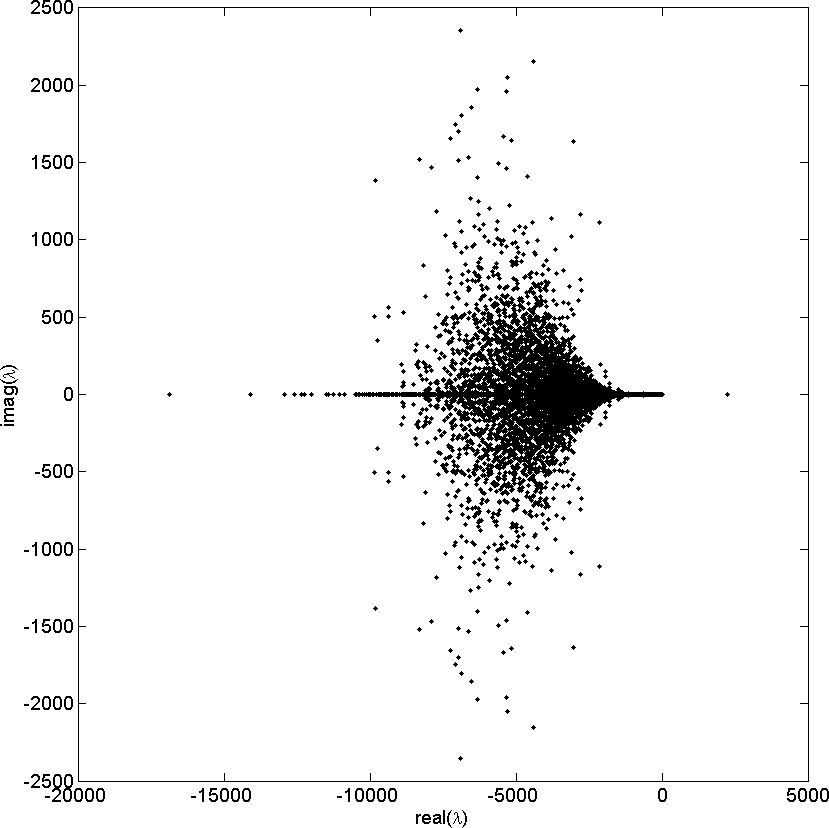}
\caption{The figure on the left of the top row shows the eigenvalues of the $n=31$ RBF-FD matrix $L_X$ for the surface Laplacian on the Red Blood Cell using $N = 10000$ MD nodes mapped to the Red Blood Cell and the per-stencil shape parameter optimization strategy with $\kappa_T=10^{12}$.  The right figure on the top row is similar, but shows the eigenvalues of $L_X$ using a single shape parameter of $\ep = 2.51$, which is the mean of the shape parameters from $L_X$ in the left figure.  The figures on the bottom row are similar to the top, but show the eigenvalues of $L_X$ for the double tours using $N = 5041$ scattered nodes.  In the left one, the per-stencil shape parameter optimization strategy was used, while the one on the right used the mean of the shape parameters from the right which was $\ep = 2.47$. }
\label{fig:eps_compare}	
\end{figure}

The idea is to choose a shape parameter $\ep_k > 0$ for each stencil $P_k$ that ``induces'' a particular target condition number $\kappa_T$ for the RBF interpolation matrix on that stencil. In the previous section we denoted this matrix by $A_{P_k}$, but now we denote it by $A_{P_k}(\ep)$ since the entries of the matrix depend continuously on the shape parameter (see Equation \eqref{eq:rbf_linsys}). The condition number of RBF interpolation matrices increase monotonically as the shape parameter decreases to zero (\emph{cf}.~\cite{FZ07}), so that the unique $\ep_k$ that induces the desired condition number $\kappa_T$ is given as the zero of the function
\begin{align}
F(\ep,\kappa_T) = \log( \kappa(A_{P_k}(\ep))/\kappa_T),
\label{eq:target_cond}
\end{align}  
where $\kappa(A_{P_k})$ is the condition number of $A_{P_k}(\ep)$ with respect to the two-norm.  Since $A_{P_k}(\ep)$ is symmetric, this is just the ratio of its largest singular value to its smallest.  We view this process as a homogenization that compensates for irregularities in the node distribution. It is a generalization of the method from~\cite{FlyerLehto2012} for the surface of the sphere, where the nodes $X$ are quasi-uniformly distributed so that one shape parameter gives roughly equal condition numbers amongst all the stencil interpolation matrices.  In that study, the shape parameter is chosen to be proportional to $\sqrt{N}$, which keeps all the conditions number approximately equal as $N$ grows.

We illustrate the effect of the proposed optimization process on the eigenvalues of the matrix approximation to the surface Laplacian $L_X$ with two tests: one on a slightly distorted but somewhat regular set of nodes and one on a very irregular set of nodes.  

For the first test, we start with the $N=10000$ quasi-uniform Maximal Determinant (MD) node set for the unit sphere (obtained from~\cite{WomerSloan:2003}). We then map this point set to an idealized Red Blood Cell surface, which is biconcave in shape; see~\cite[Appendix B]{FuselierWright2012} for the analytical expression and the upper left picture in Figure \ref{fig:rcm} for a plot of these mapped nodes.  While the MD points offer a quasi-uniform sampling of the sphere, they do not offer a good sampling when mapped to the Red Blood Cell (for a true quasi-uniform sampling of the latter, the correct procedure would be to solve an optimization problem and directly obtain MD points on the Red Blood Cell).  Next, we form two RBF-FD matrix approximations to the surface Laplacian on the Red Blood Cell using $n=31$ point stencils.  The first approximation uses an optimized shape parameter on each stencil with the target condition number set to $\kappa_T = 10^{12}$ in Equation \eqref{eq:target_cond}. The second approximation uses a single shape parameter of $\ep =  2.51$ across all stencils.  This value is the mean of the shape parameters obtained in the first approximation.  The eigenvalues of the corresponding differentiation matrices for these two procedures are shown in the top row of Figure \ref{fig:eps_compare}, with the optimized $\ep$ per stencil on the left and the single $\ep$ on the right.  We can see from the figure that optimized version produces eigenvalues all in the left half plane, while the single-$\ep$ version results in one large positive eigenvalue.

For the second test, we start with an $N=5041$ set of nodes on the double-torus that were obtained from the program 3D-XplorMath. This node set offers a fairly irregular sampling of the double-torus. For more on how the nodes were generated, see \cite{FuselierWright2012}. As on the Red Blood Cell, we form two RBF-FD matrix approximations to the surface Laplacian on the double-torus using $n=31$ point stencils. The first approximation uses an optimized shape parameter on each stencil with the target condition number set to $\kappa_T = 10^{11}$ in Equation \eqref{eq:target_cond}, the largest condition number that we could safely use on the irregular node set for $n = 31$ nodes (with larger condition numbers giving us eigenvalues with positive real parts). The second approximation uses a single shape parameter of $\ep = 2.47$ across all stencils. Again, this value is the mean of the shape parameters obtained in the first approximation. The eigenvalues of the corresponding differentiation matrices are shown in the bottom row of Figure \ref{fig:eps_compare}, with the optimized $\ep$ per stencil on the left and the single $\ep$ on the right. Again, we can see from the figure that the optimized version produces eigenvalues all in the left half plane, while the single-$\ep$ version results in one large positive eigenvalue.

We note that it is possible to choose a single shape parameter in the above examples that is sufficiently large so that $L_X$ have all eigenvalues in the left half-plane. However, this does not produce equally good results.  The reason is that smaller shape parameters generally give better accuracy (\emph{cf.}~\cite{Wright200699,LLHF}). Using the optimization procedure for selecting $\ep$ allows us to benefit from the accuracy afforded by smaller shape parameters where possible, as well as the stability afforded by larger shape parameters (when required by the irregularity of the node set). The trade-off in this procedure is that optimizing the shape parameter adds the cost of root-finding to the RBF-FD method. Additionally, fixing a target condition number across all stencils could mean that we end up choosing a lower condition number on some stencil than the condition number naturally dictated by the minimum width on that stencil. In this scenario, we are sacrificing some degree of local accuracy for the overall stability of the method. However, our tests did not reveal any impact of this on the convergence of the method.

While the shape parameter optimization procedure adds to the cost of our method, there are a few mitigating factors. First, we only solve the optimization problem to an absolute tolerance of $10^{-4}$; this proved sufficient for the purpose of stability and achieving the target condition number. Second (and more important), since the optimization is done on a per-stencil basis and the stencil computations themselves are easily parallelized, the overall optimization procedure itself is also embarassingly parallel. These advantages are retained even if the surface sampled by the node set is evolving in time.

We conclude by noting that algorithms for the stable computation of RBF-FD matrices for all value of the shape parameters are available~\cite{LLHF}, but efforts to make these work for the case when the nodes are distributed on a lower dimensional surface embedded in $\mathbb{R}^d$ is still needed. Though successful outcomes in this area would mean that we could use larger target condition numbers in our method, this does not necessarily imply the obsolescence of our optimization procedure. Given that current methods to stably compute RBF interpolants in planar domains are currently at least $5-10$ times as costly as the standard RBF interpolation method, it is probable that new methods will have this drawback as well. In such a scenario, our shape parameter optimization procedure will likely allow for cost-efficient implementations of the RBF-FD method, allowing trade-offs between accuracy and computational cost.

\section{Convergence studies}
\label{sec:convergence}
We now present the results illustrating the convergence of our method for the (forced) diffusion equation given in Equation \eqref{eq:diffusion_scalar} on some standard surfaces. We present experiments with two different optimization strategies. First, we present results for studies where we fix the condition number across all stencils for a given $N$, but allow the target condition number to grow with increasing $N$. Then, we present results for studies involving fixing the condition number for increasing $N$ (equivalent to increasing the shape parameter for increasing $N$). In the latter case, we run into saturation errors~\cite{Fasshauer:2007} due to the employment of stationary interpolation.  We use the Backward Difference Formula of Order 4 (BDF4) for all tests and set the time-step to $\Delta t = 10^{-4}$, a time-step that allows the spatial errors to dominate the temporal error. We use $\bicgstab$ to solve the implicit system arising from the BDF4 discretization; we noticed that the solver needed at most three iterations per time-step to converge to a relative tolerance of $10^{-12}$. For convenience, we measure errors using the $\ell_2$ and $\ell_{\infty}$ norms rather than approximations to the continuous versions of these norms on the test surfaces. The convergence of our method is a function of the fill distance $h_X$, defined as the radius of the largest ball that is completely contained on the manifold which does not contain a node in $X$. For quasi-uniformly distributed nodes on our test surfaces, we expect that $h \propto \frac{1}{\sqrt{N}}$, where $N$ is the total number of nodes on the surface. In the following subsections, we therefore examine convergence as a function of $\sqrt{N}$.

\subsection{Convergence studies with increasing condition number}
\label{sec:results1a}

In this section, we present the results of numerical convergence studies where the uniform condition number across the RBF-FD matrices is allowed to grow as the number of points ($N$) on the surface increases. In the absence of the shape parameter optimization procedure, this would be equivalent to fixing the shape parameter while increasing $N$, the simplest approach to take with RBF interpolation.

First, we examine the convergence of our MOL formulation by approximating the diffusion equation on a sphere. Then, we examine the convergence of our method on simulating two forced diffusion equations, one on the sphere and one on a torus. We present results for different stencil sizes $n$ and examine convergence as the total number of nodes $N$ increases. 

\vspace{0.2cm}

\textbf{Diffusion on the sphere:}

\begin{figure}[ht]	
\centering
\includegraphics[width=2.3in,height=2.2in]{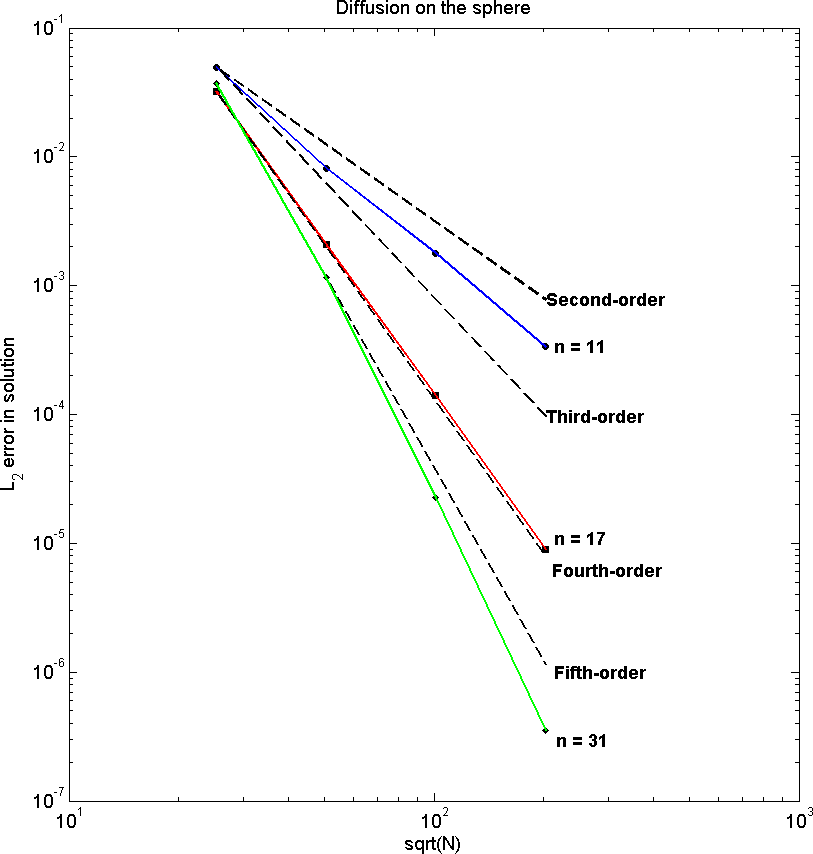}
\includegraphics[width=2.3in,height=2.2in]{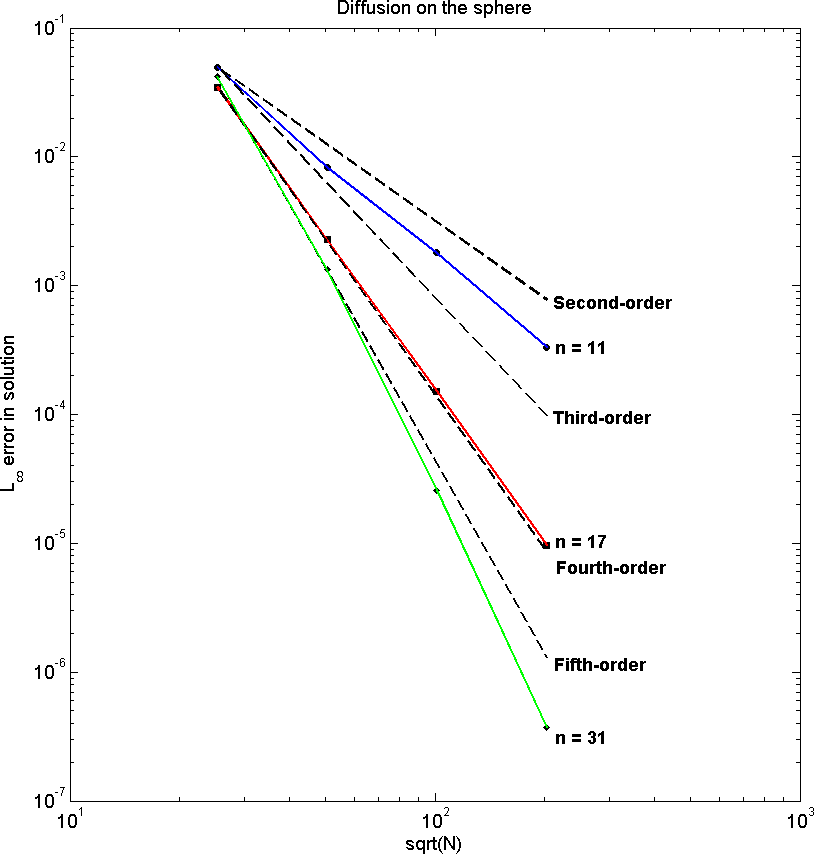}
\caption{The figure on the left shows the $\ell_2$ error in the numerical solution to the diffusion equation on the sphere as a function of $\sqrt{N}$, while the figure on the right shows the $\ell_{\infty}$ error. Both figures use a log-log scale, with colored lines indicating the errors in our method and dashed lines showing ideal $p$-order convergence for $p = 2,3,4,5$. All errors were measured against the exact solution. The errors for $n = 31$ and $N = 40962$ were computed in quad-precision.}
\label{fig:sph}	
\end{figure}

This test problem was presented in \cite{MacDonaldRuuth2009}, and involves solving the heat equation on a unit sphere $\mathbb{S}^2$. The exact solution to this problem is given as a series of spherical harmonics $u(t,\theta,\phi) = \frac{20}{3\pi}\sum_{l=1}^{\infty} e^{-l^2/9}e^{-t l (l+1)}Y_{ll}(\theta,\phi)$, where $\theta$ and $\phi$ are longitude and latitude respectively, and $Y_{lm}$ is the degree $l$ order $m$ real spherical harmonic. Since the coefficients decay rapidly, the series is truncated after $30$ terms. As in \cite{MacDonaldRuuth2009}, we evolve the PDE until $t = 0.5$, using the exact solutions to boot-strap our BDF4 scheme. We test the RBF-FD method for $n = 11,17,$ and $31$ for $N = 642, 2562, 10242,$ and $40962$ icosahedral points on the sphere~\cite{icos}, and plot the relative error in the numerical solutions. For all tests, we start with a target condition number of $\kappa_T = 10^5$ and allow it to grow with increasing $N$ to $\kappa_T = 10^{18}$. This effectively fixes the mean shape parameter on the surface as $N$ grows. The results of this study are shown in Figure \ref{fig:sph}.

In addition to the errors, Figure \ref{fig:sph} shows dashed lines corresponding to ideal $p$-order convergence, where $p = 2,3,4,5$. It is clear that our method gives convergence between orders two and three for $n=11$, close to order four for $n = 17$ and slightly higher than order five for $n=31$, both in the $\ell_2$ and the $\ell_{\infty}$ norms. Our method achieves similar results for smaller values of $N$ than were used by the Closest Point method in \cite{MacDonaldRuuth2009}, as is to be expected from a method that uses points only in the embedded space $\mathbb{S}^2$. However, it is important to be cautious when comparing errors against the Closest Point method. The values of $N$ given in the results in \cite{MacDonaldRuuth2009} are greater than the actual number of points used in that work to compute approximations to the Laplace-Beltrami operator. The values of $N$ in that work correspond to all the points used in the embedding space $\mathbb{R}^3$.

We note that the RBF-FD weights for $n = 31$ and $N = 40962$ were computed in quad-precision, though the simulations that used the weights were only run in double-precision. This is because our approach of allowing the condition number to grow with $N$ leads to condition numbers of $10^{18}$ for very high $N$ and $n$, which correspond to nearly-singular or singular matrices in double-precision. A possible way of remedying this is to start with a smaller target condition number for $N = 642$. Of course, this will lead to a higher error for each $N$, but can help offset the ill-conditioning for very large $N$ and $n$. Later in this section, we will present an alternative way of ameliorating this issue.

\vspace{0.2cm}

\textbf{Forced diffusion on the sphere:}

\begin{figure}[ht]
\centering
\includegraphics[width=2.3in,height=2.2in]{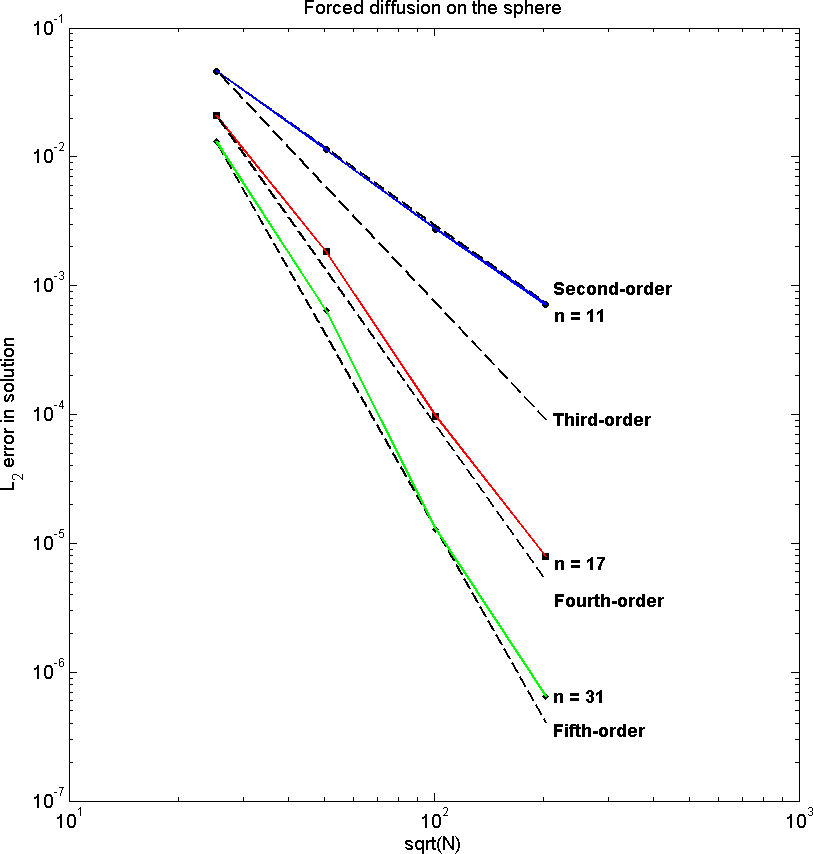}
\includegraphics[width=2.3in,height=2.2in]{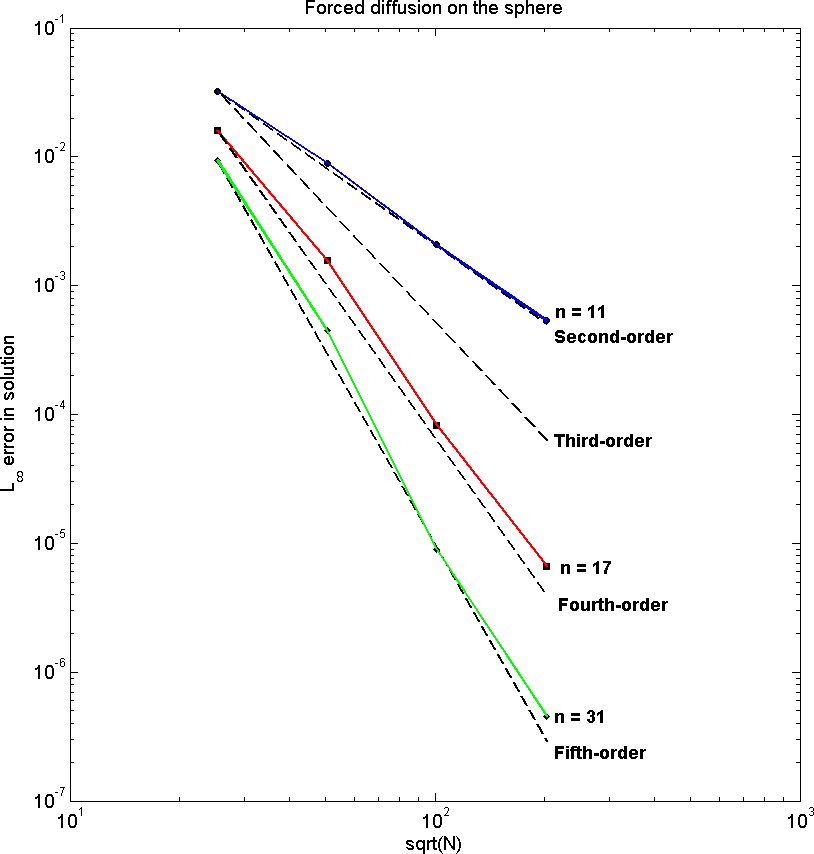}
\caption{The figure on the left shows the $\ell_2$ error in the numerical solution to the forced diffusion equation on the sphere given by Equation \eqref{eq:sph_forced} as a function of $\sqrt{N}$, while the figure on the right shows the $\ell_{\infty}$ error. Both figures use a log-log scale, with colored lines indicating the errors in our method and dashed lines showing ideal $p$-order convergence for $p = 2,3,4,5$. All errors were measured against the solution given by Equation \eqref{eq:sph_forced}. The errors for $n = 31$ and $N = 40962$ were computed in quad-precision.}
\label{fig:fsph}	
\end{figure}

This problem was first presented in \cite{Calhoun2009}, and used in \cite{FuselierWright2012} as well. For this test, we manufacture a solution to the diffusion equation on the sphere, with the forcing term $f(t,u)$ in Equation \eqref{eq:diffusion_scalar} chosen so that this solution is maintain for all time. The manufactured solution is given by
\begin{align}
u(t,\vX) = e^{-5t} \sum_{k=1}^{23} e^{-10 \cos^{-1}( {\bf \xi}_k \cdot \vX)},
\label{eq:sph_forced}
\end{align}
where ${\bf \xi}_k, k = 1,\ldots,23$ are randomly placed points on the surface of the sphere. The solution is $C^{\infty}(\mathbb{S}^2)$. As in \cite{FuselierWright2012}, we compute the forcing function analytically and evalute it implicitly in time. We compare errors in the numerical solution of the forced diffusion equation at time $t = 0.2$ for different values of $N$ and $n$. Again, we use $N = 642, 2562, 10242,$ and $40962$ icoshedral points on the sphere, with $n = 11, 17,$ and $31$ points in each of the $N$ stencils on the surface. The random placement of Gaussian centers makes this a more difficult test than diffusion of a spherical harmonic. The results are shown in Figure \ref{fig:fsph}. Again, the figure shows dashed lines corresponding to ideal $p$-order convergence, where $p = 2,3,4,5$. For this test, our method gives convergence of order two for $n=11$, close to order four for $n = 17$ and close to order five for $n=31$, both in the $\ell_2$ and the $\ell_{\infty}$ norms, which are similar to the previous experiment. Again, the weights for $n = 31$ and $N = 40962$ were computed in quad-precision, for the same reasons as before.

\vspace{0.2cm}

\textbf{Forced diffusion on a torus:}

\begin{figure}[ht]
\centering
\includegraphics[width=2.2in,height=2.2in]{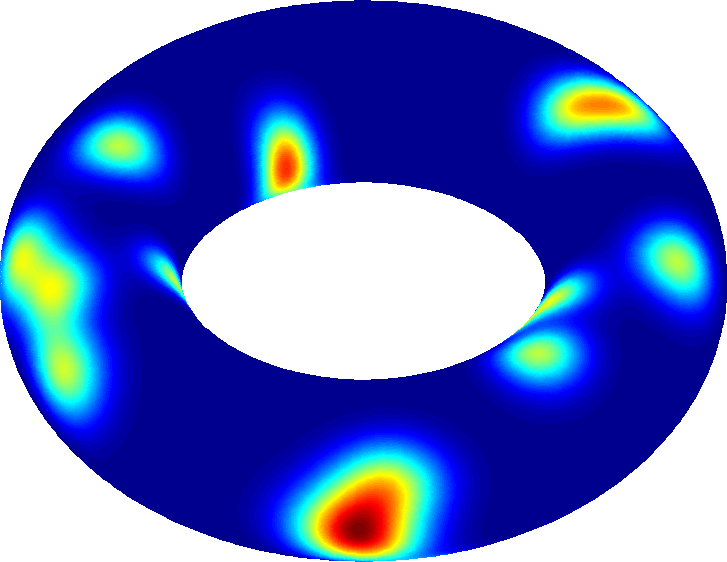}
\hspace{0.2cm}
\includegraphics[width=2.2in,height=2.2in]{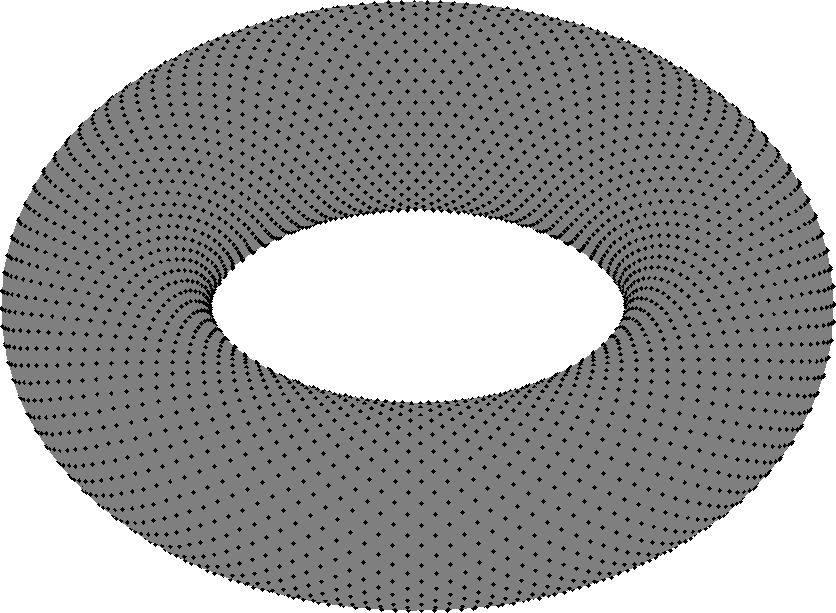}
\caption{Forced diffusion on a torus. The figure on the left shows the intial condition for the forced diffusion problem on the torus given by Equation \eqref{eq:torus}. The figure on the right shows a node set containing $N=5400$ quasi-uniformly
spaced nodes on the same torus.}
\label{fig:torpic}	
\end{figure}

This test is similar to the test involving randomly placed Gaussians on the sphere, except that this procedure is done on a torus.  We consider the torus given by the implicit equation:
\begin{align*}
\mathbb{T}^2 = \left\{\vX = (x,y,z) \in \mathbb{R}^3 \;\left| \left(1 - \sqrt{x^2 + y^2}\right)^2 + z^2 - \frac{1}{9} = 0 \right. \right\},
\end{align*}
which can be parameterized using intrinsic coordinates $\varphi$ and $\lambda$ as follows:
\begin{align}
x = \lp 1 + \frac{1}{3}\cos(\varphi)\rp\cos(\lambda),\;  y = \lp 1 + \frac{1}{3}\cos(\varphi)\rp\sin(\lambda),\; z = \frac{1}{3}\sin(\varphi),
\label{eq:torus}
\end{align}
where $-\pi\leq \varphi,\lambda\leq \pi$.
The surface Laplacian of a scalar function $f:\mathbb{T}\rightarrow\mathbb{R}$ in this intrinsic coordinate system is given as
\begin{align*}
\Delta_{\mathbb{M}} f(\varphi,\lambda) = \frac{1}{\lp 1 + \frac{1}{3}\cos(\varphi)\rp^2} \frac{\partial^2 f}{\partial \lambda^2} + \frac{9}{\lp 1 + \frac{1}{3}\cos(\varphi)\rp}  \frac{\partial}{\partial \varphi}\left(\lp 1 + \frac{1}{3}\cos(\varphi)\rp\frac{\partial f}{\partial \varphi}\right).
\end{align*}
The manufactured solution to the diffusion equation (given by Equation \eqref{eq:diffusion_scalar} with $\mathbb{M} = \mathbb{T}$) is
\begin{align}
u(t,\varphi,\lambda) = e^{-5t} \sum_{k=1}^{23} e^{-a^2(1-\cos(\lambda-\lambda_k))-b^2(1-\cos(\varphi-\varphi_k))},
\label{eq:torus_forced}
\end{align}
where $a=9$, $b=3$, and $(\varphi_k,\lambda_k)$ are randomly chosen values in $[-\pi,\pi]^2$. The solution is $C^{\infty}(\mathbb{T}^2)$ and a visualization at $t=0$ is given in Figure \ref{fig:torpic} (left).  While the solution and forcing function are all specified using intrinsic coordinates, the RBF-FD method uses only extrinsic (Cartesian coordinates) without requiring knowledge of the underlying intrinsic coordinate system. As with the forced diffusion test on the sphere, we compute the forcing function corresponding to Equation \eqref{eq:torus_forced} analytically and evaluate it implicitly. We similarly compare errors in the numerical solution of the forced diffusion equation at time $t = 0.2$ for stencils of size $n = 11, 17$ and $31$ nodes.
\begin{figure}[ht]
\centering
\includegraphics[width=2.3in,height=2.2in]{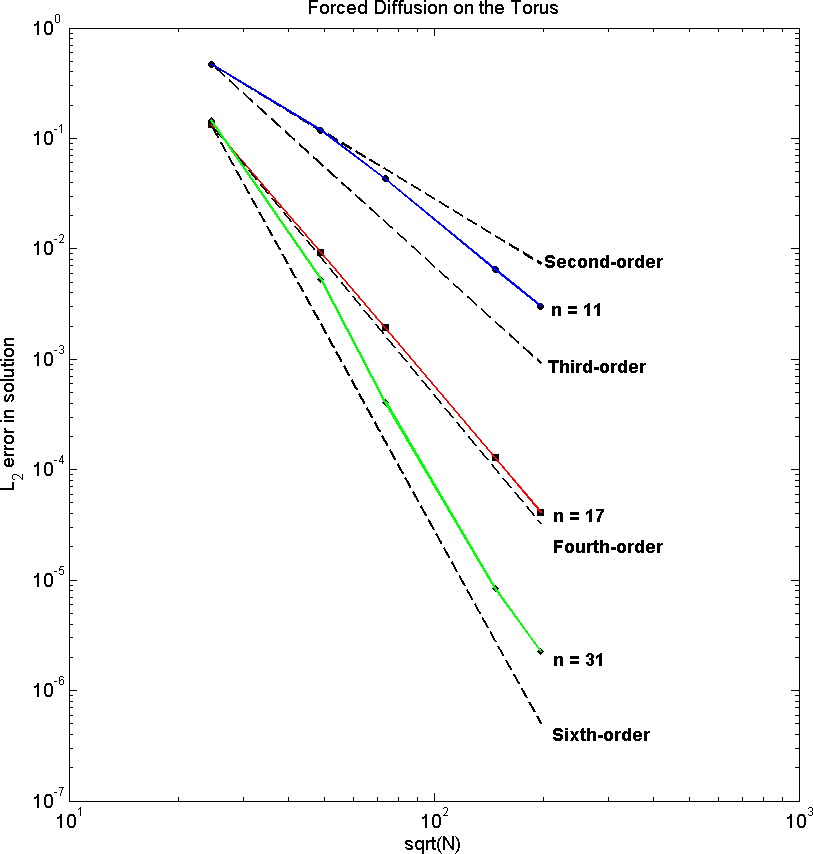}
\includegraphics[width=2.3in,height=2.2in]{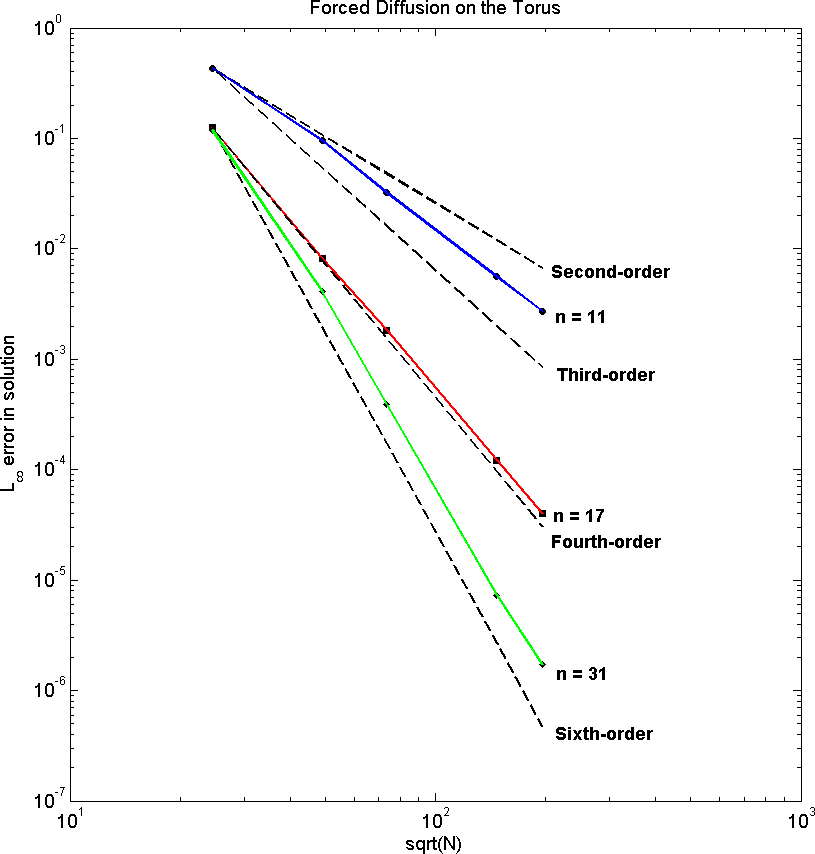}
\caption{The figure on the left shows the $\ell_2$ error in the numerical solution to the forced diffusion equation on the torus given by Equation \eqref{eq:torus} as a function of $\sqrt{N}$, while the figure on the right shows the $\ell_{\infty}$ error. Both figures use a log-log scale, with colored lines indicating the errors in our method and dashed lines showing ideal $p$-order convergence for $p = 2,3,4,6$. All errors were measured against the solution given by Equation \eqref{eq:torus_forced}. The errors for $n = 31$ and $N = 38400$ were computed in quad-precision.}
\label{fig:ftor1}	
\end{figure}
The node sets we use for experiments on the torus are generated from a ``staggered'' grid in intrinsic variable space, and are determined as follows:
\begin{enumerate}
\item Given $m$, choose $m+1$ equally spaced angles on $[-\pi,\pi]$ in $\varphi$ and $3m+1$ equally spaced angles on $[-\pi,\pi]$ in $\lambda$.
\item Disregard the values of $\varphi$ and $\lambda$ at $\pi$ and take a direct product of the remaining points to obtain $N = 3m^2$ points on $[-\pi,\pi)^2$. Map these points to $\mathbb{T}$ using Equation \eqref{eq:torus} and call the set of nodes $X_1$.
\item Next, generate another set of $N=3m^2$ gridded points in $[-\pi,\pi)^2$ from the previous set by offsetting the $\varphi$ coordinate by $\pi/m$ and the $\lambda$ coordinate by $\pi/(3m)$, so they lie at the midpoints of the previous gridded values. Map these to $\mathbb{T}$ and call the set of nodes $X_2$.
\item The final set of nodes is given by $X=X_1 \cup X_2$.  
\end{enumerate}
In the experiments we use $m=10, 20, 30, 60, 80$, corresponding to node sets of size $N = 600, 2400, 5400, 21600, 38400$.  A plot of the nodes for $N=5400$ is shown in Figure \ref{fig:torpic} (right).  These points remain more or less uniformly spaced on the torus as $N$ grows.

The results for the experiments are shown in Figure \ref{fig:ftor1}. Again, the figure shows dashed lines corresponding to ideal $p$-order convergence, where $p = 2,3,4,6$. On this test, our method gives convergence of order two for $n=11$, close to order four for $n = 17$ and between orders five and six for $n=31$, both in the $\ell_2$ and the $\ell_{\infty}$ norms. The convergence rates are comparable to the results seen for diffusion on the sphere, and slightly better than those seen for the forced diffusion problem on the sphere. The results for $n = 31$ and $N = 38400$ were computed in quad-precision, for the same reasons as before.

\subsection{Convergence studies with fixed condition number}
\label{sec:results1b}

In Section \ref{sec:results1a}, we saw that allowing the condition number to grow as $N$ increases by fixing the mean shape parameter can give excellent results, but will eventually cause the RBF interpolation matrices to be ill-conditioning for large $N$ and $n$.

In this section, we present an alternate approach. We choose to fix the target condition number at a particular (reasonably large) value for all values of $N$ and $n$. As $N$ increases, this has the effect of increasing the value of the average shape parameter. Our goal here is to understand the relationship between the magnitude of the target condition number and the value of $n$ and $N$ at which saturation errors can set in. This would also give us intuition on the connection between the target condition number and the order of convergence of our method.

With this in mind, we present the results of numerical convergence studies of forced diffusion on a sphere and on a torus for two fixed condition numbers, $\kappa_T = 10^{14}$ and $\kappa_T = 10^{20}$, for increasing values of $N$ and $n$. For $\kappa_T = 10^{20}$, the RBF-FD weights on each patch were run in quad-precision using the Advanpix Multicomputing Toolbox. However, once the weights were obtained, they were converted back to double-precision and the simulations were carried out only in double-precision.

\begin{figure}[ht]
\centering
\includegraphics[width=2.3in,height=2.2in]{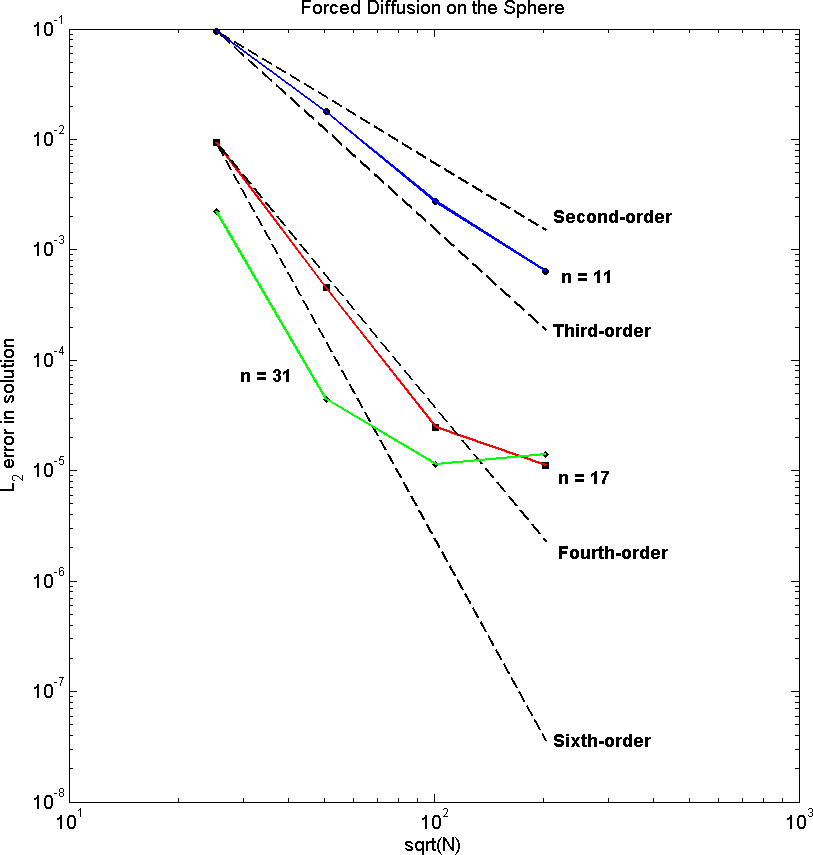}
\includegraphics[width=2.3in,height=2.2in]{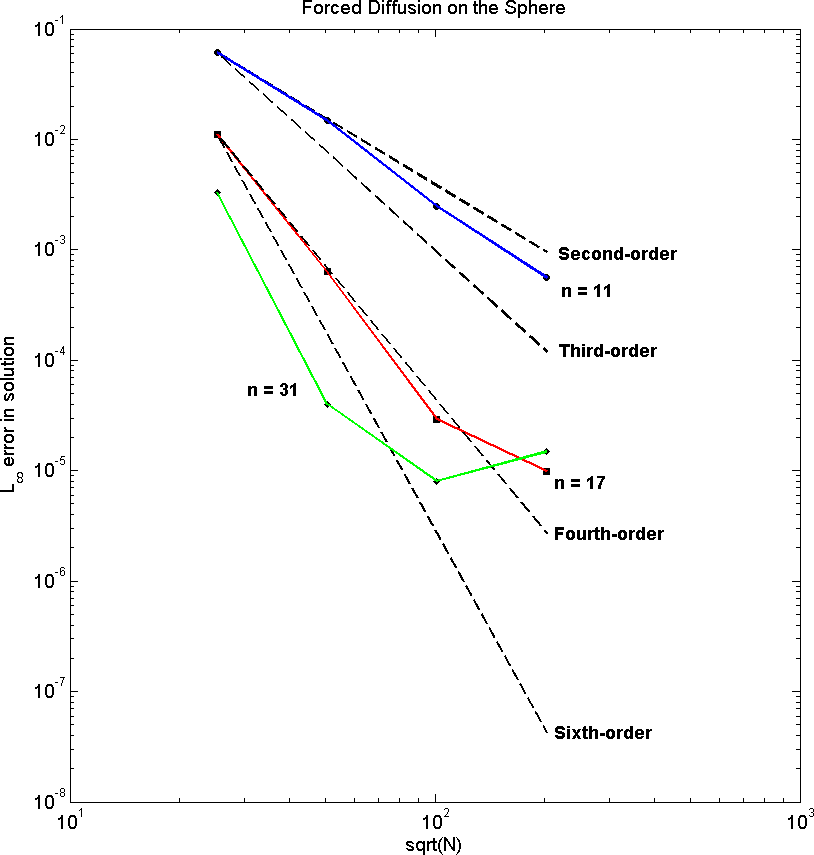}
\caption{The figure on the left shows the $\ell_2$ error in the numerical solution to the forced diffusion equation on the sphere given by Equation \eqref{eq:sph_forced} as a function of $\sqrt{N}$, while the figure on the right shows the $\ell_{\infty}$ error. Both figures use a log-log scale, with colored lines indicating the errors in our method and dashed lines showing ideal $p$-order convergence for $p = 2,3,4,6$. The target condition number was set to $\kappa_T = 10^{14}$.}
\label{fig:fsphc2}	
\end{figure}

\begin{figure}[ht]
\centering
\includegraphics[width=2.3in,height=2.2in]{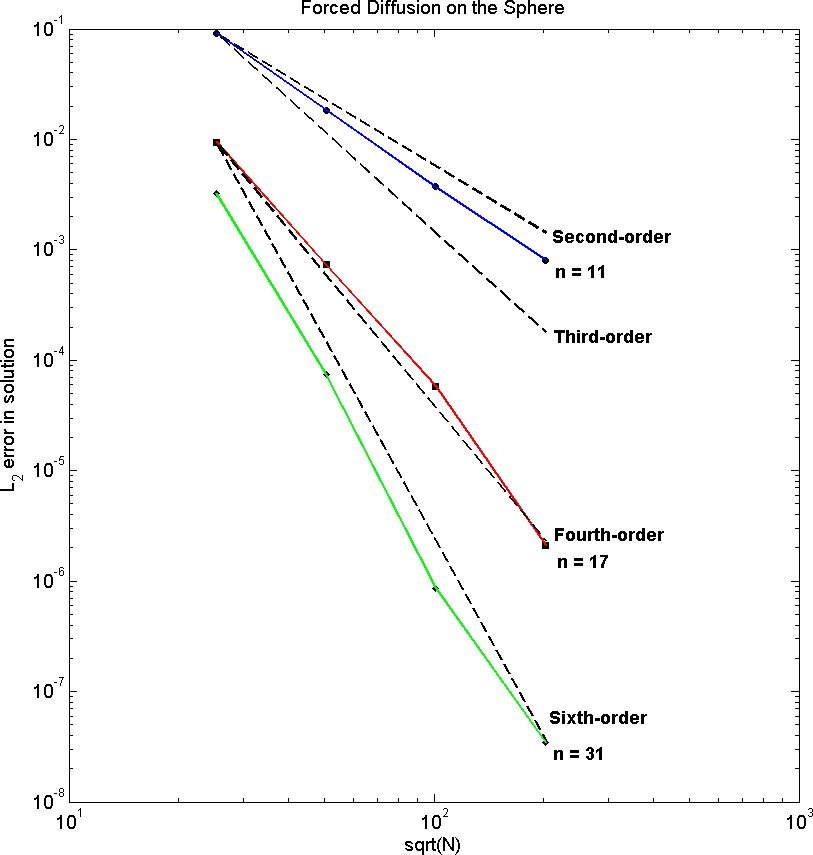}
\includegraphics[width=2.3in,height=2.2in]{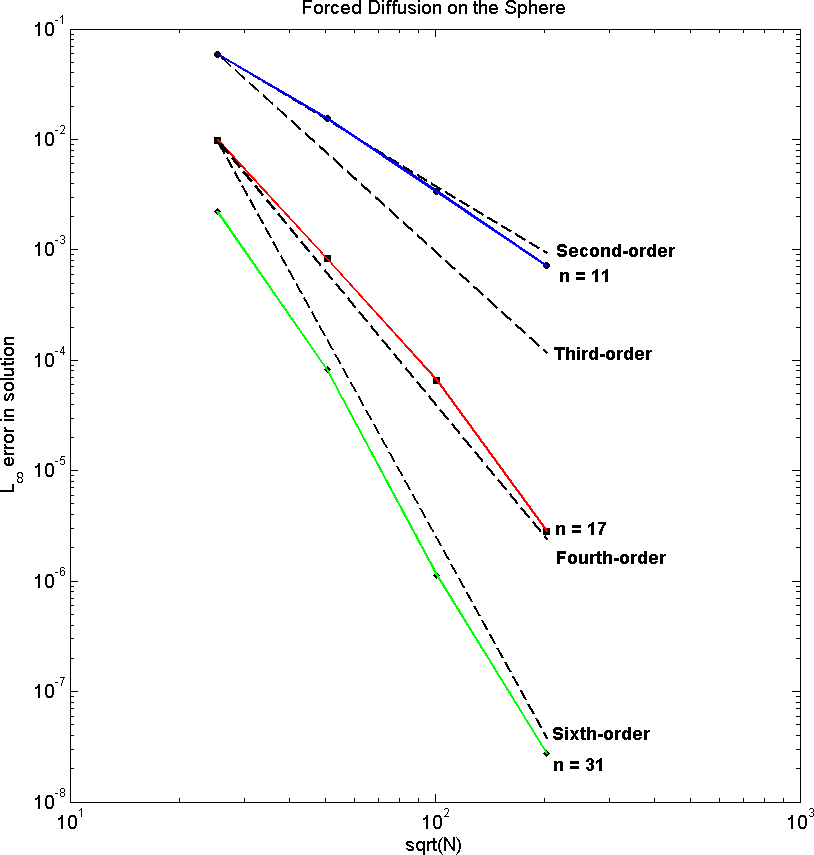}
\caption{The figure on the left shows the $\ell_2$ error in the numerical solution to the forced diffusion equation on the sphere given by Equation \eqref{eq:sph_forced} as a function of $\sqrt{N}$, while the figure on the right shows the $\ell_{\infty}$ error.  Both figures use a log-log scale, with colored lines indicating the errors in our method and dashed lines showing ideal $p$-order convergence for $p = 2,3,4,6$. The target condition number was set to $\kappa_T = 10^{20}$ and all RBF-FD weights were computed in quad precision.}
\label{fig:fsphc3}	
\end{figure}

\begin{figure}[ht]
\centering
\includegraphics[width=2.3in,height=2.2in]{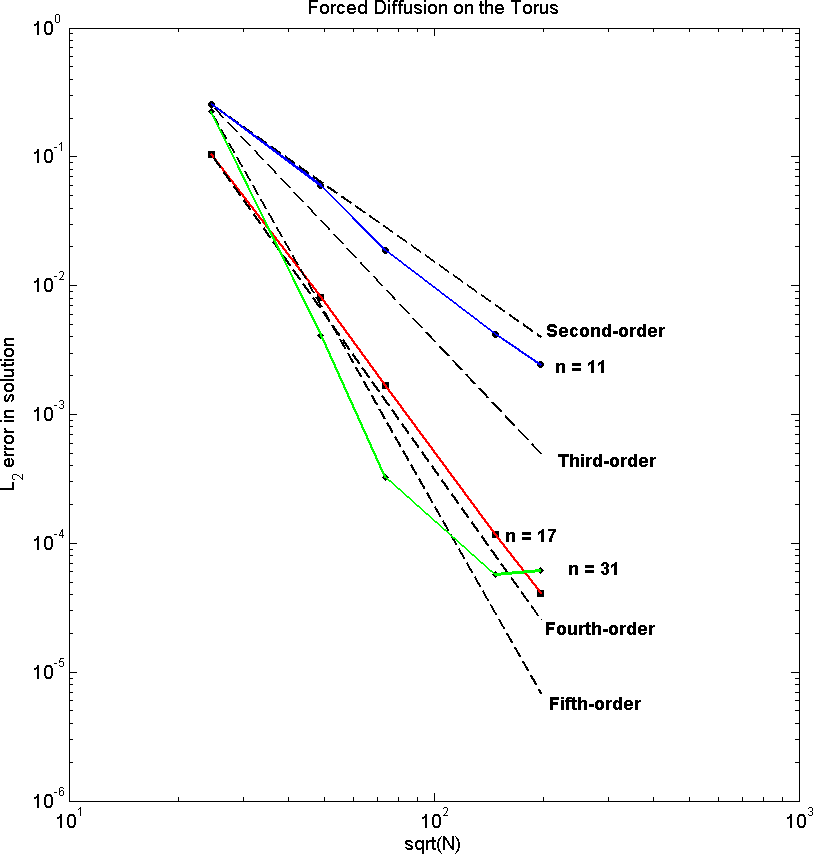}
\includegraphics[width=2.3in,height=2.2in]{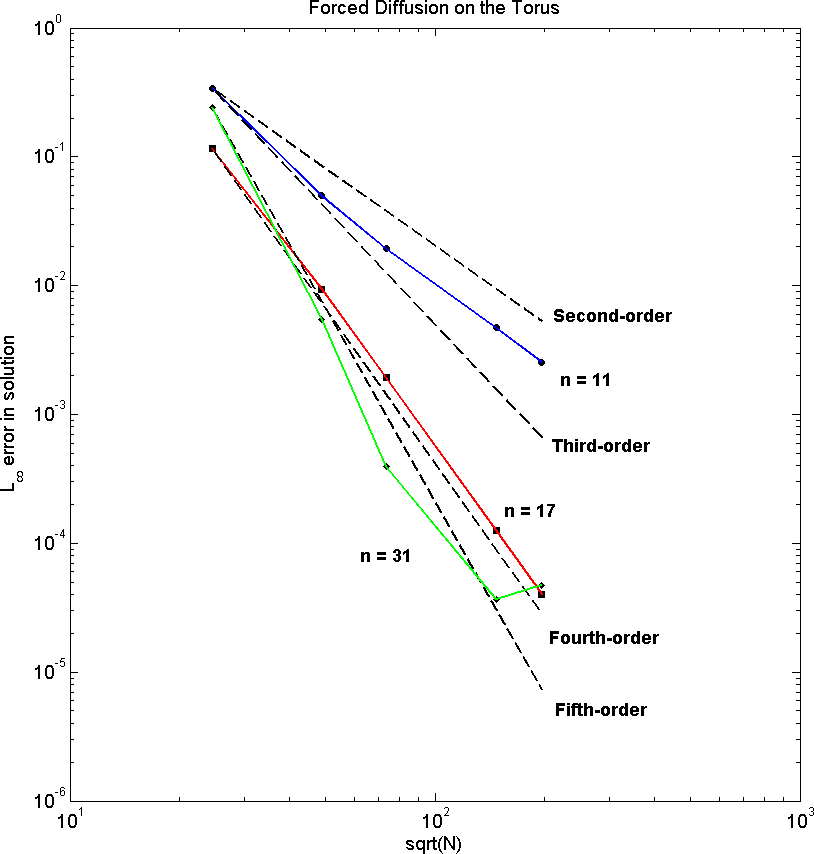}
\caption{The figure on the left shows the $\ell_2$ error in the numerical solution to the forced diffusion equation on the torus given by Equation \eqref{eq:torus} as a function of $\sqrt{N}$, while the figure on the right shows the $\ell_{\infty}$ error. Both figures use a log-log scale, with colored lines indicating the errors in our method and dashed lines showing ideal $p$-order convergence for $p = 2,3,4,5$. The target condition number was set to $\kappa_T = 10^{14}$.}
\label{fig:ftorc2}	
\end{figure}

\begin{figure}[ht]
\centering
\includegraphics[width=2.3in,height=2.2in]{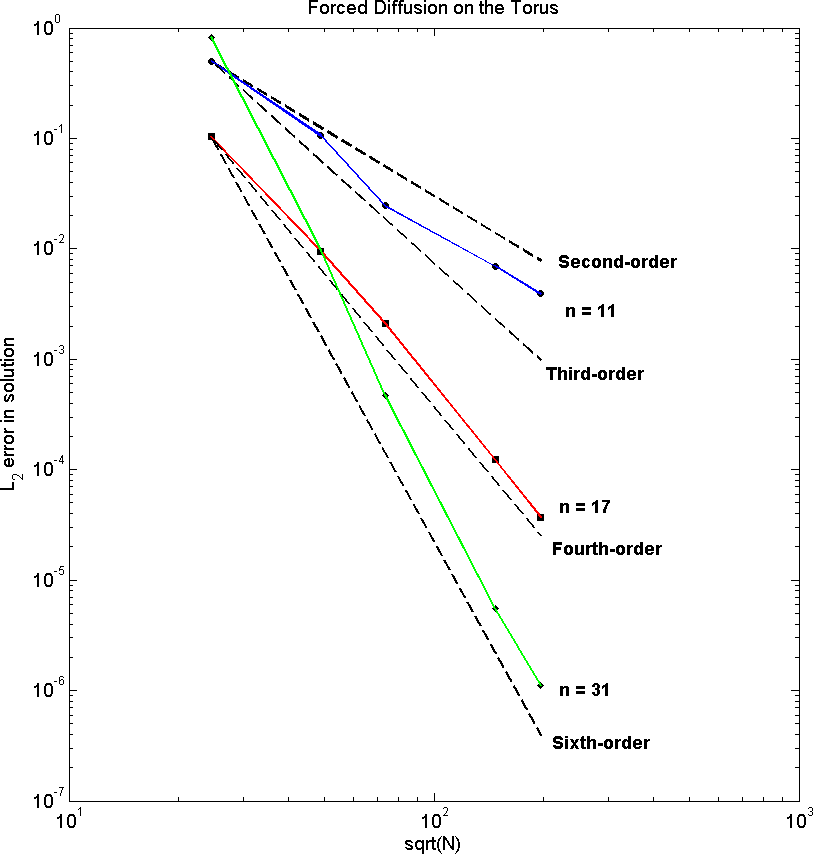}
\includegraphics[width=2.3in,height=2.2in]{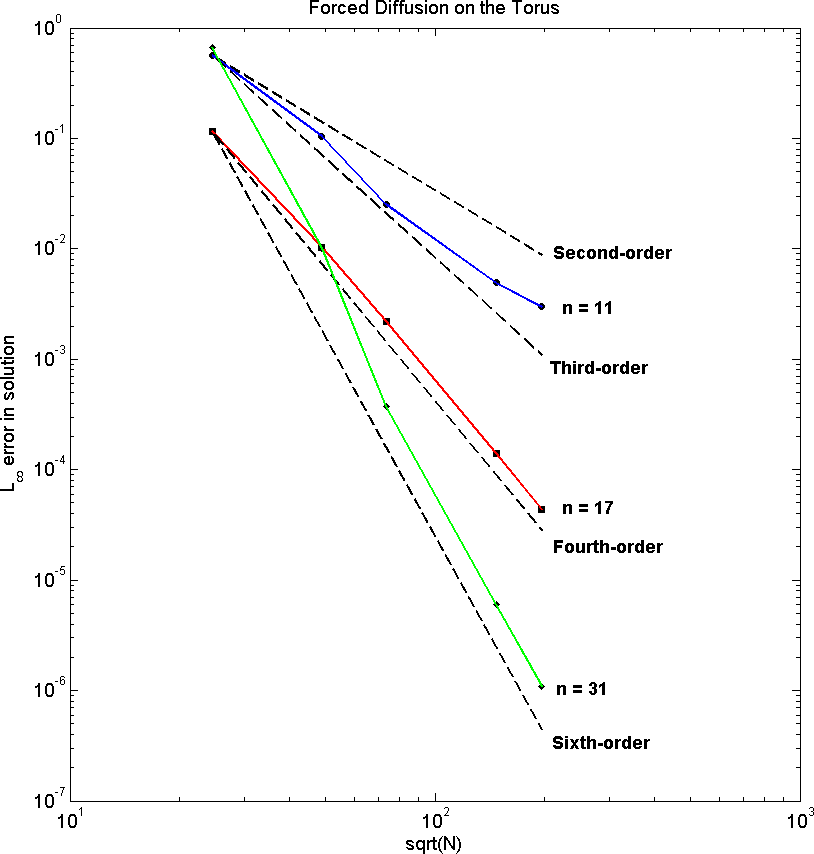}
\caption{The figure on the left shows the $\ell_2$ error in the numerical solution to the forced diffusion equation on the torus given by Equation \eqref{eq:torus} as a function of $\sqrt{N}$, while the figure on the right shows the $\ell_{\infty}$ error. Both figures use a log-log scale, with colored lines indicating the errors in our method and dashed lines showing ideal $p$-order convergence for $p = 2,3,4,6$. The target condition number was set to $\kappa_T = 10^{20}$ and all RBF-FD weights were computed in quad precision.}
\label{fig:ftorc3}	
\end{figure}

We test the behavior of our MOL formulation for forced diffusion on the sphere and for forced diffusion on the torus, with the solution to the former given by Equation \eqref{eq:sph_forced} and the solution to the latter given by Equation \eqref{eq:torus_forced}. We again use a BDF4 method with $\Delta t = 10^{-4}$, with the forcing term computed analytically and evaluated implicitly. The errors in the $\ell_2$ and $\ell_{\infty}$ measured at $t = 0.2$ are shown in Figures \ref{fig:fsphc2}--\ref{fig:fsphc3} for the sphere, and in Figures \ref{fig:ftorc2}--\ref{fig:ftorc3} for the torus.

\vspace{0.2cm}

\textbf{Forced diffusion on a sphere:}

The results for $\kappa_T = 10^{14}$ are shown in Figure \ref{fig:fsphc2}, and those for $\kappa_T = 10^{20}$ are shown in Figure \ref{fig:fsphc3}. Figure \ref{fig:fsphc2} shows that fixing the target condition number at $\kappa_T = 10^{14}$ produces no saturation errors in the $\ell_2$ or $\ell_{\infty}$ norms for $n = 11$. In fact, this value of $\kappa_T$ appears to produce saturation in the $\ell_2$ norm only for large values of $N$ for $n = 17$. We see similar saturation errors in the $\ell_{\infty}$ norm. For $n = 31$, we see saturation for $N > 2562$; this is a clear indication that generating high-order RBF-FD methods requires the ability to use large target condition numbers (small values of the shape parameter $\ep$).

It is clear from Figure \ref{fig:fsphc3} that there are no saturation errors for $n = 11$ and $n = 17$ for the values of $N$ tested when $\kappa_T = 10^{20}$. In addition, for $n = 31$, we see convergence that is close to order six with no saturation errors for the values of $N$ used. This confirms that using small shape parameters within the RBF-FD method can give higher convergence. This is further motivation for the development of algorithms to stably compute the RBF interpolation matrices as $\ep \to 0$.

\vspace{0.2cm}

\textbf{Forced diffusion on a torus:}

The results for $\kappa_T = 10^{14}$ are shown in Figure \ref{fig:ftorc2}, and those for $\kappa_T = 10^{20}$ are shown in Figure \ref{fig:ftorc3}. Figure \ref{fig:ftorc2} shows that fixing the target condition number at $\kappa_T = 10^{14}$ produces no saturation errors in the $\ell_2$ or $\ell_{\infty}$ norms for $n = 11$ or $n = 17$ in either norm, in contrast to forced diffusion on the sphere. Indeed, $\kappa_T = 10^{14}$ seems sufficient for methods up to order 4 for the values of $N$ tested. However, for $n = 31$, we see saturation errors for $N > 5400$, again showing that larger target condition numbers are required for high-order RBF-FD methods.

Figure \ref{fig:ftorc3} shows that convergence on the torus is a bit more erratic than on the sphere for the forced diffusion problem when $\kappa_T = 10^{20}$. While there are no saturation errors for $n = 11$ and $n = 17$ for the values of $N$ tested, the convergence is slightly lower for smaller values of $N$, possibly indicating that the node set on the torus is not quite as uniformly-spaced as the one on the sphere for those values of $N$. However, as $N$ is increased, the rate of convergence seems to be slightly better than what was seen on the sphere. This becomes apparent when looking at the line for $n = 31$. We see large errors for $N = 600$, but a rapid fall-off as $N$ is increased, with the overall order of convergence for $n = 31$ being between five and six.

\section{Application: Turing patterns}
\label{sec:results2}

This section presents an application of our RBF-FD method to solving a two-species Turing system (two coupled reaction-diffusion equations) on different surfaces. We present two types of results, the first is for surfaces where parameterizations or implicit equations describing the surface are known, and the second where they are not.
\begin{table}[ht]
\centering
\begin{tabular}{|c|c|c|c|c|c|c|c|}
\hline
\textbf{Surface/Pattern} & $\delta_v$ & $\alpha$ & $\beta$ & $\gamma$ & $\tau_1$ & $\tau_2$ & Final time \\ \toprule \hline
RBC/spots & $4.5 \times 10^{-3}$ & $0.899$ & $-0.91$ & $-0.899$ & $0.02$ & $0.2$ & 800 \\\hline
RBC/stripes & $2.1 \times 10^{-3}$ & $0.899$ & $-0.91$ & $-0.899$ & $3.5$ & $0$ & 6500\\\hline
Bumpy sphere/spots & $4.5 \times 10^{-3}$ & $0.899$ & $-0.91$ & $-0.899$ & $0.02$ & $0.2$ & 800 \\\hline
Bumpy sphere/stripes & $2.1 \times 10^{-3}$ & $0.899$ & $-0.91$ & $-0.899$ & $3.5$ & $0$ & 7000 \\\hline
Double-torus/spots & $2.1 \times 10^{-3}$ & $0.899$ & $-0.91$ & $-0.899$ & $0.02$ & $0.2$ & 700 \\\hline
Double-torus/stripes & $8.87 \times 10^{-4}$ & $0.899$ & $-0.91$ & $-0.899$ & $3.5$ & $0$ & 6000\\\hline
Frog/spots & $2.87 \times 10^{-4}$ & $0.899$ & $-0.91$ & $-0.899$ & $0.02$ & $0.2$ & 600\\\hline
Bunny/stripes & $2.87 \times 10^{-4}$& $0.899$ & $-0.91$ & $-0.899$ & $3.5$ & $0$ & 6000 \\\hline
\end{tabular}
\caption{The above table shows the values of the parameters of Equations \eqref{eq:turing1} and \eqref{eq:turing2} used in the numerical experiments shown in Figures \ref{fig:tur_spots_1} and \ref{fig:tur_spots_2}. In all cases, we set $\delta_u = 0.516 \delta_v$.}
\label{tab:tur_param}
\end{table}
To facilitate comparison, we use the Turing system first described for the surface of the sphere in~\cite{VareaEtAl1999} and applied to more general surfaces in~\cite{FuselierWright2012}. The system describes the interaction of an activator $u$ and inhibitor $v$ according to
\begin{align}
\frac{\partial u}{\partial t} &=  \alpha u(1- \tau_1 v^2) + v(1 - \tau_2 u) + \delta_u\Delta_{\mathbb{M}} u, \label{eq:turing1} \\ 
\frac{\partial v}{\partial t} &= \beta v \lf(1 + \frac{\alpha \tau_1}{\beta} uv \rt) + u(\gamma + \tau_2v) + \delta_v \Delta_{\mathbb{M}} v. \label{eq:turing2}
\end{align}
If $\alpha = -\gamma$, then $(u,v) = (0,0)$ is a unique equilibrium point of this system. Altering the diffusivity rates of $u$ and $v$ can lead to instabilities which manifest as pattern formations. The coupling parameter $\tau_1$ favors stripe formations, while $\tau_2$ favors spots. Stripe formations take much longer to attain ``steady-state'' than spot formations. In the following subsections,  We use the Semi-implicit Backward Difference Formula of order 2 (SBDF2) as the time-stepping scheme, and set the time-step to $\Delta t = 0.01$ for all tests. Since the diffusion terms are handled implicitly, the RBF-FD matrix needs to be inverted every time-step. We accomplish this by pre-computing a sparse LU decomposition of the matrix, and using the triangular factors for forward and back solves every time-step.  The values for all parameters for Equations \eqref{eq:turing1} and \eqref{eq:turing2}, including final times for simulations, are presented in Table \ref{tab:tur_param}.
\begin{table}[ht]
\centering
\begin{tabular}{|c|c|c|c|}
\hline
\textbf{Surface} & Number of nodes ($N$)  & Stencil size ($n$) & Target Cond. No. ($\kappa_T$) \\ \toprule \hline
RBC & $10000$ & $31$ & $10^{12}$\\\hline
Bumpy sphere & $10000$ & $31$ & $10^{12}$ \\\hline
Double-torus & $12100$ & $31$ & $10^{11}$ \\\hline
Frog & $7458$ & $31$ & $10^{10}$ \\\hline
Bunny & $11339$ & $31$ &$ 10^{10}$ \\\hline
\end{tabular}
\caption{The above table shows the values of the parameters used in the RBF-FD discretization of Equations \eqref{eq:turing1} and \eqref{eq:turing2} for the numerical experiments shown in Figures \ref{fig:tur_spots_1} and \ref{fig:tur_spots_2}. In all cases, the time-step was set to $\Delta t = 0.01$.}
\label{tab:tur_simparam}
\end{table}

\subsection{Turing patterns on manifolds}

\begin{figure}[htbp]
\begin{center}
        \includegraphics[scale = 0.3]{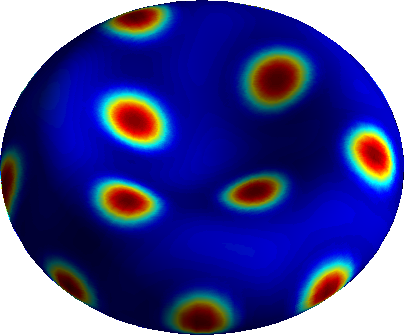}
        \hspace{1cm}
        \includegraphics[scale = 0.33]{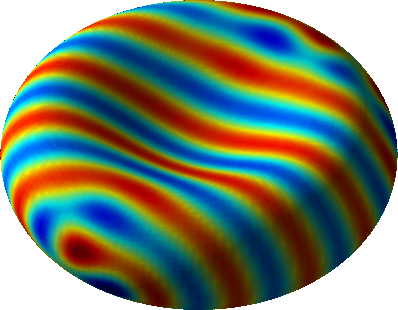}  
               				
				\vspace{0.5cm}        				
        \includegraphics[scale = 0.25]{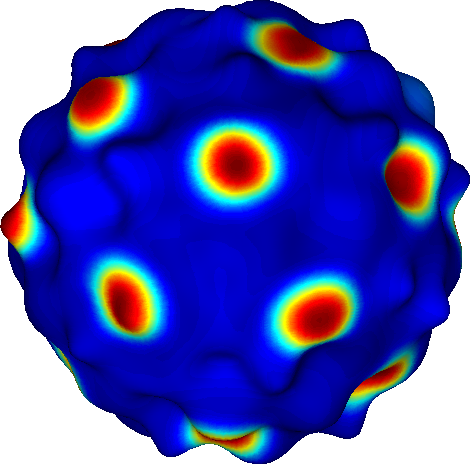}
        \hspace{1cm}
        \includegraphics[scale = 0.25]{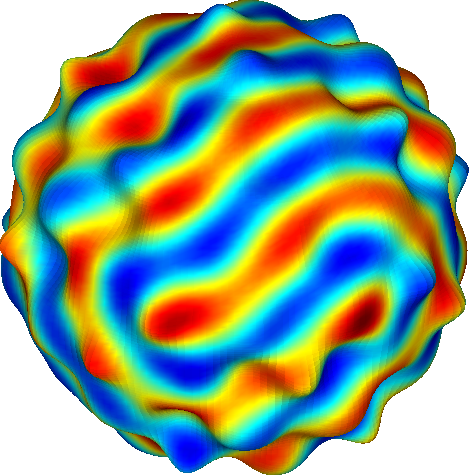}
        
				\vspace{0.5cm} 				
        \includegraphics[scale = 0.19]{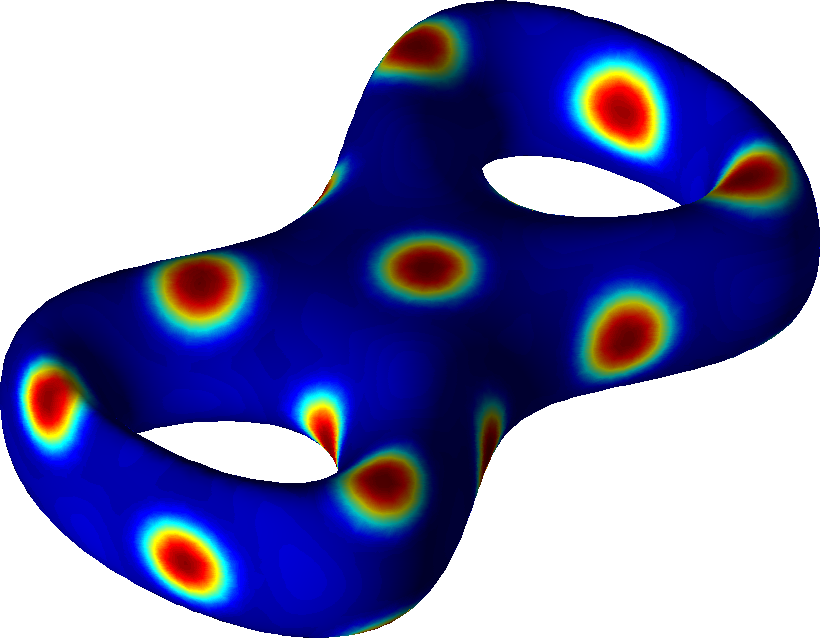}
        \hspace{1cm}
        \includegraphics[scale = 0.18]{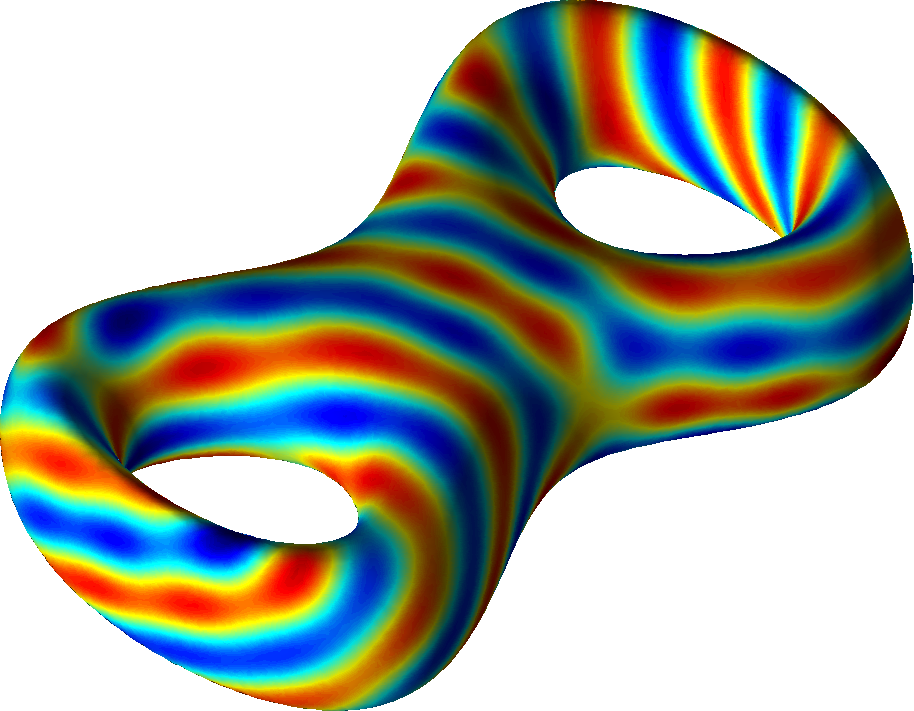}
        
        \caption{Steady Turing spots and stripe patterns resulting from solving Equations \eqref{eq:turing1} and \eqref{eq:turing2} on the Red Blood Cell, the Bumpy Sphere model and the double-torus surfaces. In all plots, red corresponds to a high concentration of $u$ and blue to a low concentration.}
        \label{fig:tur_spots_1}
\end{center}
\end{figure}

We first solve the Turing system on three surfaces: the Red Blood Cell (RBC) and the double-torus described earlier, and on the Bumpy Sphere detailed in \cite{FuselierWright2012}. RBCs are biconcave surfaces and can be represented parametrically, as described earlier. The Bumpy Sphere is a point set downloaded from an online repository, and equipped with point unit normals by parametric interpolation with the RBF parametric model presented in \cite{Shankar2013}. Also, since the downloaded model had $N = 5256$ vertices and we wished to demonstrate the viability of the RBF-FD method on far more vertices, we sampled our parametric model to generate $N = 10000$ vertices, and solve the Turing system on that point set. The first three rows of Table \ref{tab:tur_simparam} list all the parameters used in the RBF-FD discretizations of Equations \eqref{eq:turing1} and \eqref{eq:turing2} on each of these three surfaces, and the last column of Table \ref{tab:tur_param} lists the final times used for these simulations. The results of these simulations are shown in Figure \ref{fig:tur_spots_1}. The spot and stripe patterns are qualitatively similar to those shown in \cite{FuselierWright2012}.

\subsection{Turing patterns on more general surfaces}

\begin{figure}[ht]
\begin{center}
        \includegraphics[scale = 0.33]{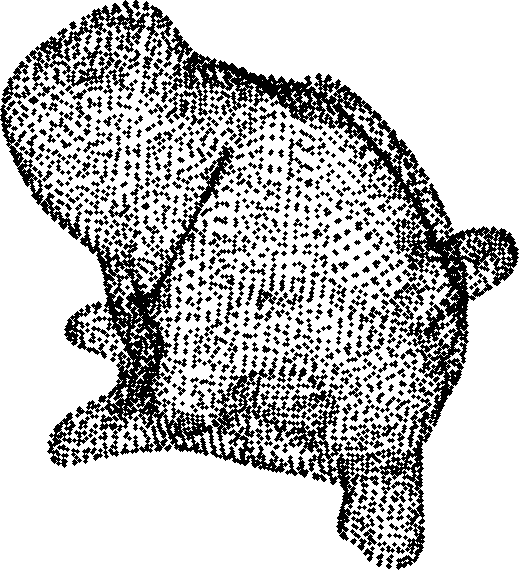}
        \hspace{1cm}
        \includegraphics[scale = 0.4]{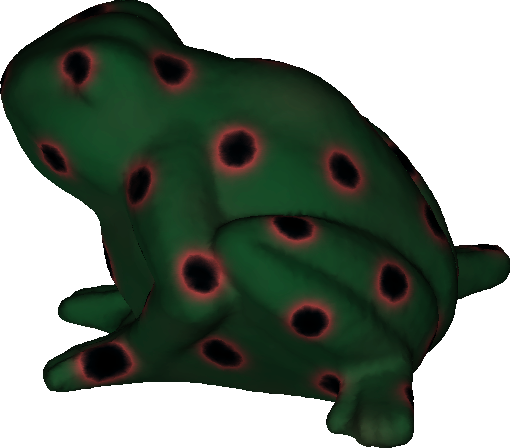}        
        \vspace{0.5cm}
        \includegraphics[scale = 0.25]{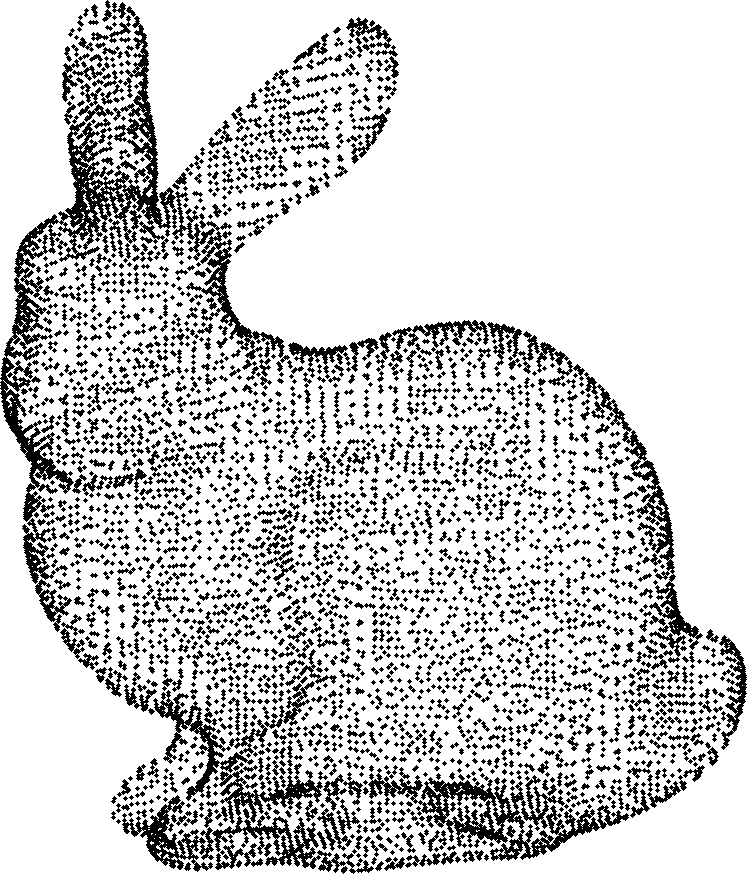}  
        \hspace{1cm}
        \includegraphics[scale = 0.25]{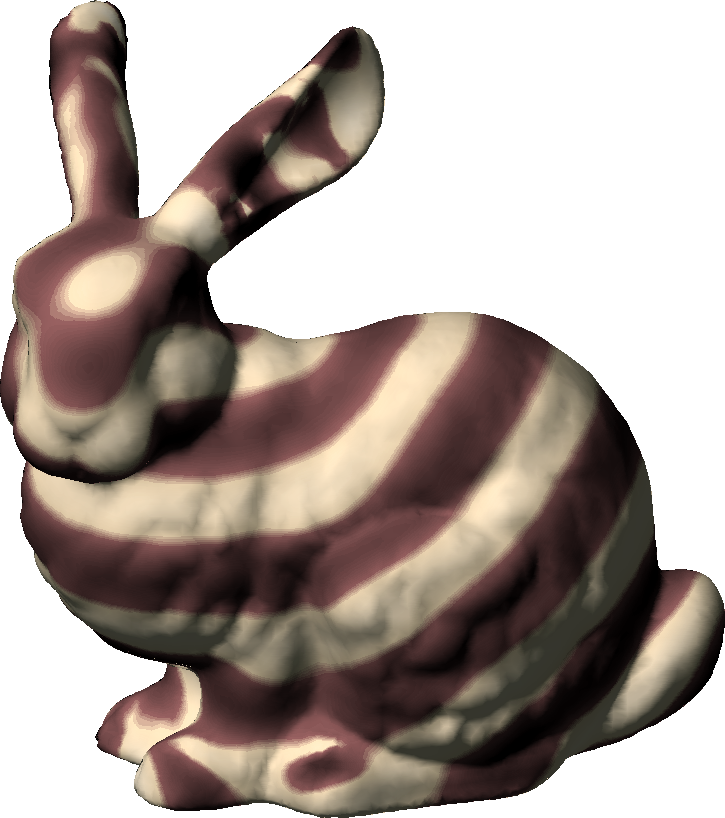}  
        \caption{The figure on the top right shows a Turing spot pattern on a Frog model. Green corresponds to a low concentration, and brown and black to higher concentrations. The figure on the bottom right shows a Turing stripe pattern on the Stanford Bunny model. Here, the lightest browns (almost white) correspond to low concentrations, and darker browns correspond to higher concentrations. Both the figures on the left show the point clouds used for the solution of the Turing system.}
        \label{fig:tur_spots_2}
\end{center}
\end{figure}

We now turn our attention to more general point sets: a Frog model (obtained from the AIM@SHAPE Shape Repository) and the Stanford Bunny model (obtained from the Stanford 3D Scanning Repository). Rather than as point clouds, these models are available in the form of meshes and approximate normal data. In contrast the previous two examples, it is not clear if the surfaces represented by these meshes can be analytically parametrized. To prepare these point sets for simulations, we first run the Poisson surface reconstruction algorithm~\cite{Kazhdan:2006} to generate a water-tight implicit surface that fits the point cloud. This algorithm requires both the point cloud and the approximate normals as input. Forming an implicit surface smoothes the approximate normal vectors input into the Poisson surface reconstruction, resulting in a more stable RBF-FD discretization. Having generated an implicit surface, we sample that with the Poisson disk sampling algorithm to generate a point cloud with the desired number of points. The other rationale for employing Poisson disk sampling is that while the Poisson surface reconstruction will fix any holes in the mesh, those former holes may not be sufficiently sampled. This pre-processing was performed entirely in MeshLab~\cite{CCR08}. 

After this preprocessing, we run a Turing spot simulation on the Frog model, and a Turing stripe simulation on the Stanford Bunny. The last two rows of Table \ref{tab:tur_simparam} list all the parameters used in the RBF-FD discretizations of Equations \eqref{eq:turing1} and \eqref{eq:turing2} on each of these surfaces, and the last column of Table \ref{tab:tur_param} lists the final times used for these simulations. The results are shown in Figure \ref{fig:tur_spots_2}. Before the color-mapping for aesthetics, the results are qualitatively similar to those shown in Figure \ref{fig:tur_spots_1}.

\section{Discussion}
\label{sec:discuss}

In this paper, we introduced a new numerical method based on Radial Basis Function-generated Finite Differences (RBF-FD) for computing a discrete approximation to the Laplace-Beltrami operator on surfaces of codimension one embedded in $\mathbb{R}^3$. The method uses scattered nodes on the surface, without requiring expansion into the embedding space. We improved on the method presented in \cite{ShankarWrightEtAlIJNMF2013}, designing a stable numerical method that does not require stabilization with artificial viscosity (a feature of RBF-FD methods for convective flows). This development was facilitated by an algorithm to optimize the shape parameter for each interpolation patch on the surface. We demonstrated that this optimization procedure can compensate for irregularities in the sampling of the surface. We then presented error and convergence estimates for our method using two approaches: allowing the condition number to grow with the number of points on the surface, and fixing the condition number for an increasing number of points. We discussed the trade-offs inherent in each approach, and provided intuition as to the relationship between the condition number, shape parameter and the order of convergence of our method on the diffusion equation on a sphere and a torus. We presented an application of our method to simulating reaction-diffusion equations on surfaces; specifically, we demonstrated the solution of Turing PDEs on several interesting shapes, both parametrizable and more general.

While our method currently works for static objects, our goal is to apply RBF-FD to the solution of PDEs on evolving surfaces, with the evolution dictated by the interaction of a fluid with the object. It will be necessary to employ efficient and fast k-d tree implementations, including algorithms for dynamically updating and/or re-balancing the k-d tree as the point set evolves. The method will almost certainly need to be parallelized to be efficient.

One issue with our method is its ability to handle thin features on surfaces. Indeed, our method is not innately robust to such features. For RBF-FD to be robust on more general surfaces, it will be necessary to combine our method with an adaptive refinement code that detects thin features and samples sides of the feature sufficiently (the alternative would be to find an efficient way to measure distances along arbitrary surfaces).

A natural extension of this work would be to adapt the method to handle spatially-variable (possibly anisotropic) diffusion. While this extension is not conceptually difficult, the realization of this extension would make our method even more useful for biological applications, like the simulation of gels or viscoelastic materials on surfaces. We intend to address this in a follow-up study. Finally, while we have successfully applied RBF-FD to periodic surfaces, it would be interesting to apply the method to solving PDEs on surfaces with boundary conditions imposed on them. We intend to address this in a follow-up study as well.

\section*{Acknowledgments}
We would like to acknowledge useful discussions concerning this work 
within the CLOT group at the University of Utah. The first, third and
fourth authors acknowledge funding support under NIGMS grant R01-GM090203.  
The second author acknowledges funding support under NSF-DMS grant 
1160379 and NSF-DMS grant 0934581.
 
\bibliographystyle{spmpsci}
\bibliography{article}

\end{document}